\documentclass{amsart}
\UseRawInputEncoding

\def\sg{\sigma}
\def\D{\Delta}
\def\la{\lambda}
\def\e{\varepsilon}

\usepackage{graphicx}
\usepackage{subfigure}
\usepackage{booktabs} 
\usepackage{multirow}
\usepackage{array}
\usepackage{slashed}
\usepackage{amsmath,bm}
\usepackage{mathrsfs,psfrag,eepic,epsfig}
\usepackage{fancyhdr}
\usepackage{makecell}
\usepackage{float}
\usepackage{amsmath,amssymb}
\usepackage{epstopdf}
\usepackage{caption}

\newtheorem{theorem}{Theorem}[section]
\newtheorem{lemma}{Lemma}[section]
\newtheorem{corollary}{Corollary}[section]

\newtheorem{remark}{Remark}[section]

\begin{document}

%
%

\title[Convergence of the GLTR method]{On convergence of the generalized Lanczos
trust-region method for trust-region
subproblems}
\author{Bo Feng}
\address{China University of Mining and Technology, Xuzhou, 221116, P.R. China.}
\curraddr{}
\email{bofeng@cumt.edu.cn}
\thanks{}

\author{Gang Wu}
\address{China University of Mining and Technology, Xuzhou, 221116, P.R. China.}
\curraddr{}
\email{gangwu@cumt.edu.cn}
\thanks{This
author is supported by the National Natural Science Foundation of China under grant 12271518,
the Key Research and Development Project of Xuzhou Natural Science Foundation under grant
KC22288, and the Open Project of Key Laboratory of Data Science and Intelligence Education of
the Ministry of Education under grant DSIE202203.}

\begin{abstract}
The generalized Lanczos trust-region (GLTR) method is one of the most popular approaches for solving large-scale trust-region subproblem (TRS). Recently, Z. Jia \emph{et al.} [Z. Jia and F. Wang, \emph{SIAM J. Optim., 31 (2021), pp. 887--914}] considered the convergence of this method and established some {\it a prior} error bounds on the residual, the solution and the Lagrange multiplier. In this paper, we revisit the convergence of the GLTR method and try to improve these bounds. First, we establish a sharper upper bound on the residual. Second, we give a new bound on the distance between the approximation and the exact solution, and show that the convergence of the approximation has nothing to do with the spectral separation. Third, we present some \emph{non-asymptotic} bounds for the convergence of the Lagrange multiplier, and define a factor that plays an important role on the convergence of the Lagrange multiplier. Fourth, we revisit the convergence of the Krylov subspace method for the cubic regularization variant of trust-region subproblem, and substantially improve the convergence result established in [X. Jia, X. Liang, C. Shen and L. Zhang, \emph{SIAM J. Matrix Anal. Appl. 43 (2022), pp. 812--839}] on the multiplier. Numerical experiments demonstrate the effectiveness of our theoretical results.
\end{abstract}
\pagestyle{myheadings} \thispagestyle{fancy} \markboth{B. FENG AND G. WU}{\sc }

\maketitle

\section{Introduction}
In this paper, we are interested in the following large-scale \emph{trust-region subproblem} (TRS) \cite{9,7}
\begin{equation}
\begin{aligned}\label{1}
\min_{\|\bm x\|_2\leq \Delta} {\Big\{}f(\bm x)=\frac{1}{2}\bm x^T A\bm x+\bm x^T\bm g{\Big\}},
\end{aligned}
\end{equation}
where $A\in \mathbb{R}^{n\times n}$ is a symmetric matrix, $\bm g\in \mathbb{R}^n$ is a nonzero vector, and $\Delta>0$ is the trust-region radius.

The TRS \eqref{1} arises from many problems such as Tikhonov regularization of the ill-posed problem \cite{M1,M2,M3}, graph partitioning problems \cite{W.W.1}, and the Levenberg--Marquardt algorithm for solving nonlinear least squares problems \cite{7}. Moreover, solving TRS is a key step in trust-region methods for dealing with general nonlinear optimization problems \cite{9,7}. Some regularization variations of the TRS such as the cubic regularization variant of TRS \cite{Y.N} were consider in \cite{T.B,0,110,G.C.P,C.G.M,J&Z,F.L,Y.N}, which can be applied to solve unconstrained minimization problems.

There are many methods have been proposed for solving \eqref{1}.
For instance, the Mor$\acute{\rm e}$--Sorensen \cite{28} method is an efficient approach for small and dense TRS.
For large-scale problem, Gould \emph{et al}. modified the Mor$\acute{\rm e}$--Sorensen method by using the Taylor series \cite{11}, and Steihaug \emph{et al.} \cite{T.S,P.L.} solved the TRS \eqref{1} via a truncated conjugate gradient (tCG) method.
The generalized Lanczos trust-region method (GLTR) \cite{10} solves \eqref{1} by using a projection method.
A nested Lanczos method was proposed in \cite{Z.L}, which can be viewed as an accelerated GLTR method.
A matrix-free method was presented in \cite{M1,M2}, and a semi-definite programming (SDP) based method was proposed in
\cite{24}.

In \cite{GGM}, Gander, Golub, and von Matt solved \eqref{1}
from an eigenproblem point of view. More precisely,
the TRS \eqref{1} was rewritten as a standard eigenvalue problem of size $2n$.
In \cite{3}, Adachi \emph{et al.} generalized this strategy to a weighted-norm TRS which reduces to
a generalized eigenvalue problem of size $2n$. Moreover, Adachi \emph{et al.} showed that the Lagrange multiplier corresponding to the optimal solution is the rightmost eigenvalue of a matrix pair of size $2n$-by-$2n$, and the optimal solution can be computed from the eigenvector associated with the rightmost eigenvalue.
One refers to \cite{9,J.B,W.W,7,soren,Tao} and the references therein for more methods on the TRS \eqref{1}.

Indeed, the generalized Lanczos trust-region method (GLTR) \cite{10} is one of the most popular approaches for solving large-scale TRS \eqref{1}.
The convergence of GLTR was investigated in \cite{0,110,G.S.,5,4}. In this paper, we denote by $\bm x_{opt}$ an optimal solution to TRS \eqref{1}, and by $\la_{opt}$ the Lagrangian multiplier  associated with the solution $\bm x_{opt}$. Let $\bm x_k$ and $\la_{k}$ be the approximate solutions of $\bm x_{opt}$ and $\la_{opt}$ obtained from the $k$-th step of the GLTR method, respectively. Recently, some upper bounds on $f(\bm x_k)-f(\bm x_{opt})$ and $\|(A+\la_k I)\bm x_k+\bm g\|_2$ were established in \cite{0,110} and \cite{G.S.}. In \cite{4}, Zhang \emph{et al.} presented some upper bounds on $f(\bm x_k)-f(\bm x_{opt})$ and $\|\bm x_{opt}-\bm x_k\|_2$. The results showed that the convergence of the GLTR method is similar to that of the conjugate gradient (CG) method. Recently, Z. Jia \emph{et al.} \cite{4}
established some {\it a prior} error bounds on  $\vert\la_{opt}-\la_k\vert$, $\|(A+\la_k I)\bm x_k+\bm g\|_2$, $\sin\angle(\bm x_{opt},\bm x_k)$, as well as $f(\bm x_k)-f(\bm x_{opt})$.
The convergence of the Krylov subspace method for cubic regularization variant of TRS was analyzed in \cite{0,110,G.S.,J&Z}.

In this paper, we revisit the convergence of the GLTR method, and try to refine the bounds due to Z. Jia \emph{et al.} \cite{4}. Moreover, we substantially improve some results of X. Jia {\it et al.} \cite{J&Z} and Gould {\it et al.} \cite{G.S.} on the Krylov subspace method for cubic regularization variant of TRS.

The contributions of this work are as follows:
 \begin{itemize}
 \item
First, for the residual $\|(A+\la_kI)\bm x_k+\bm g\|$, we establish a \emph{sharper} upper bound than those of Z. Jia \emph{et al.} \cite[Theorem 5.7]{5} and Golud \emph{et al.} \cite[Theorem 3.4]{G.S.}. Our bound is $\mathcal{O}(\kappa^{\frac{1}{2}})$ times smaller than that of Z. Jia \emph{et al.}, where $\kappa$ is the spectral condition number of $A+\la_{opt}I$.
\item
Second, for the approximate solution $\bm x_k$, we give a new bound on $\sin\angle(\bm x_{opt},\bm x_k)$, and show that the convergence of $\bm x_k$ has nothing to do with the spectral separation $sep(\la_{opt}, C_k)$. As $sep(\la_{opt}, C_k)$ can be very small or even close to zero in practice, our new upper bound is much better than some existing ones.
 \item
Third, for the approximate multiplier $\la_k$, the error bound on $\la_{opt}-\la_k$ is \emph{asymptotic} \cite[Theorem 4.11]{5}, i.e., it holds when \emph{k} is \emph{sufficiently large}. In this paper, we present some \emph{non-asymptotic} bounds that hold for \emph{all} $k\geq 1$.
 In particular, we show that the factor $s(\lambda_{opt})$ defined in \eqref{eq513} plays an important role on
 the convergence of $\la_k$. To the best of our knowledge, it seems that this is the first analysis from this point of view.
\item
Fourth, we consider the convergence of the Krylov subspace method for solving the cubic regularization variant \eqref{2} of the TRS \eqref{1}. We substantially improve the convergence result on the multiplier $\mu_{opt}$ due to X. Jia \emph{et al.} \cite{J&Z}.

\end{itemize}

This paper is organized as follows.
In Section 2, we present some preliminaries.
In Section 3 and Section 4, we give insight into the convergence of the GLTR method, and establish some refined upper bounds on $\|(A+\la_kI)\bm x_k+\bm g\|$, $\sin\angle(\bm x_{opt},\bm x_k)$ and $\la_{opt}-\la_k$, respectively.
In Section 5, we revisit the convergence of the Krylov subspace method for the cubic regularization variant of the
TRS \eqref{1}. Numerical experiments are performed in Section 6, which illustrate the effectiveness and sharpness of the proposed new upper bounds.
Some concluding remarks are given in Section 7.

Throughout this paper, we denote by $(\cdot)^T$ the transpose of a matrix or vector, by $\mathcal{N}(A)$ the null space of a matrix $A$,
by $\mathcal{R}(A)$ the range space of a matrix $A$,
by $\|\cdot\|$ the Euclidean norm of a
matrix or vector.
Let $\mathcal{W}$ be a linear subspace of $\mathbb{R}^n$, and let $W$ be an orthonormal basis of $\mathcal{W}$. Then the distance between a nonzero vector $\bm p$ and $\mathcal{W}$ is defined as \cite{Y.S}
 \begin{equation}\label{eq22.005}
 \sin\angle(\bm p, \mathcal{W})=\min_{\bm x\in \mathcal{W}}\left\| \frac{\bm p}{\|\bm p\|}-\bm x \right\|=\frac{\|(I-WW^T)\bm p\|}{\|\bm p\|},
 \end{equation}
and the sine of the angle between a nonzero vector $\bm p$ and $\bm q$ is defines as \cite{Y.S}
 \begin{equation*}
 \sin\angle(\bm p, \bm q)=\min_{\omega\in \mathbb{C}}\left\| \frac{\bm p}{\|\bm p\|}-\omega\bm q \right\|.
 \end{equation*}
Let $\bm 0$, $\bm O$ and $ I$ be the zero
vector, zero matrix and identity matrix, respectively, whose sizes are clear from the context. Let $\bm e_i$ be the $i$-th column of the identity matrix $I$.
In this paper, $A\succeq\bm O~(A\succ\bm O)$ implies that $A$ is symmetric semi-positive definite (positive definite), and $A\succcurlyeq B~(A\succ B)$ stands for
$A-B$ is a symmetric semi-positive definite (positive definite) matrix.


\section{Preliminaries}
In this section, we briefly introduce the solutions to the subproblem \eqref{1}, the generalized Lanczos trust-region (GLTR) method \cite{10}, and some existing convergence results on the GLTR method.

\subsection{Solutions of the subproblems}
Let $A = U\Lambda U^T$ be the eigendecomposition, where $U=\begin{pmatrix}\bm u_1,\bm u_2,\ldots,\bm u_n\end{pmatrix}\in \mathbb{R}^{n\times n}$ is orthonormal, and
$$
\Lambda=diag(\alpha_1,\alpha_2,\ldots,\alpha_n)
$$
with $\alpha_1\geq\alpha_2\geq\cdots\geq\alpha_n$ being the eigenvalues of $A$ in decreasing order.
The following result provides a necessary and sufficient condition for a global optimal solution of TRS \eqref{1}.
\begin{theorem}\cite{9,28,C.sec}\label{7.59}
The vector $\bm x_{opt}$ is a global optimal solution of the trust-region problem \eqref{1}, if and only if $\|\bm x_{opt}\|\leq \Delta$ and there exists Lagrange
multiplier $\lambda_{opt} \geq 0$ such that
 \begin{equation}\label{eq2010}
(A+\lambda_{opt}I){\bm x_{opt}}=-\bm g,~\lambda_{opt}(\Delta-\|\bm x_{opt}\|)=0~~{\rm and}~~
A+\lambda_{opt}I \succcurlyeq \bm O.
\end{equation}
\end{theorem}

In particular, there are two situations for the TRS \eqref{1} \cite{W.W,28,7}:

${\bf Case ~1}.$  $\emph{Easy case}$:
$\la_{opt}>-\alpha_n$ and $\la_{opt}\geq 0$.
In this case, $A+\la_{opt}I\succ\bm O$, the solution $\bm x_{opt}$ for TRS \eqref{1} is unique and
 $
 \bm x_{opt} = -(A+\la_{opt}I)^{-1}\bm g.
 $

 ${\bf Case~2 }.$ $\emph{ Hard case}$: $\la_{opt}=-\alpha_n$. In this case, we have
 $$
 \bm g~\bot~\mathcal{N}(A-\alpha_nI)~~{\rm and}~~\|(A-\alpha_nI)^\dag\bm g\|\leq \D,
 $$
where the superscript $\dag$ denotes the Moore-Penrose inverse of a matrix.

 \subsection{The generalized Lanczos trust-region (GLTR) method}

As $A$ is a symmetric matrix, one can use the $k$-step Lanczos process to generate an orthonormal basis $Q_k=[\bm q_1,\bm q_2,\ldots,\bm q_{k}]$ for the Krylov subspace $\mathcal{K}_k(A,\bm g)=span \{ \bm g, A\bm g,\ldots,A^{k-1}\bm g\}$. Moreover, the following Lanczos relation holds \cite{2,Y.S}
 \begin{align}\label{21.18}
 AQ_k&=Q_kT_k+\beta_{k}\cdot\bm q_{k+1}\bm e_k^T,
 \end{align}
 where $\bm e_k$ is the $k$-th column of the identity matrix, and
 $$
 T_k =Q^T_kAQ_k=\begin{pmatrix}
 \delta_1 & \beta_1 & &\\
 \beta_1  &  \delta_2& \beta_2 & \\
          &  \ddots &  \ddots& \ddots&    \\
          &         & \beta_{k-2}& \delta_{k-1} & \beta_{k-1}\\
          &         &            &  \beta_{k-1} &  \delta_{k}
 \end{pmatrix}
 \in \mathbb{R}^{k\times k},
 $$
with $\bm g=\|\bm g\|\bm q_1$ and $Q^T_k \bm g=\|\bm g\|\bm e_1$.

The GLTR method due to Gould {\it et al.} \cite{10} is one of the most popular approaches for solving the large-scale TRS problem \eqref{1}.
The GLTR method solves the following TRS problem
 \begin{equation}\label{9.00}
 \begin{aligned}
 \min_{\|\bm x\|\leq\D,~\bm x\in \mathcal{K}_k(A,\bm g)}\Big\{f(\bm x)&=\frac{1}{2}\bm x^T A\bm x+\bm x^T\bm g\Big\},
 \end{aligned}
 \end{equation}
 which leads to the following {\it reduced} TRS problem
 \begin{equation}\label{9.25}
 \begin{aligned}
   \min_{ \|\bm h\|\leq \D}\left\{f_k(\bm h)= \frac{1}{2}\bm h^T T_k\bm h+\|\bm g\|\cdot\bm h^T \bm e_1\right\}.
   \end{aligned}
 \end{equation}
 Let $\bm h_k$ be the minimizer of \eqref{9.25}. It can be verified that
 \begin{equation}\label{9.000}
 \begin{aligned}
 \bm x_k= Q_k\bm h_k=\min_{\|\bm x\|\leq\D,~\bm x\in \mathcal{K}_k(A,\bm g)}f(\bm x),
 \end{aligned}
 \end{equation}
 which can be used as an approximation to $\bm x_{opt}$.

In \cite{5,4}, it was shown that $\|\bm x_k-\bm x_{opt}\|\rightarrow 0$ as $k$ increase.
Since \eqref{9.25} is also a trust-region subproblem, by Theorem \ref{7.59}, we have that $\|\bm h_k\|\leq \D$ and there exists a Lagrange multiplier $\lambda_{k} \geq 0$, such that
 \begin{equation*}
(T_k+\la_kI)\bm h_k=-\|\bm g\|\bm e_1,~\la_k(\D-\|\bm h_k\|)=0~~{\rm and}~~T_k+\la_kI\succcurlyeq\bm O.
\end{equation*}
So we have from \eqref{21.18} that
\begin{align}
 (A+\lambda_kI)\bm x_k=&(A+\la_kI)Q_k\bm h_k=Q_k(T_k+\lambda_kI)\bm h_k+ \beta_k\bm q_{k+1}\bm e_k^T\bm h_k\nonumber\\
                      =&-\|\bm g\|\cdot Q_k\bm e_1+  (\beta_k\bm e_k^T\bm h_k)\cdot\bm q_{k+1}=-{\bm g}+(\beta_k\bm e_k^T\bm h_k)\cdot\bm q_{k+1}.\label{7.55}
\end{align}
As a result, GLTR is an orthogonal projection method satisfying
 \begin{equation}\label{11.10}
\left\{
\begin{array}{ll}
\bm x_k \in  \mathcal{K}_k(A,\bm g),                 \\
\bm r_k=(A+\la_kI)\bm x_k+\bm g~ \bot ~\mathcal{K}_k(A,\bm g),
\end{array}\right.
\end{equation}
which is a key property utilized in our analysis; see Section 3.
Notice that $T_k$ is a tridiagonal and irreducible matrix, the TRS \eqref{9.25} is also in the \emph{easy case} \cite[Theorem 5.3]{10}, and the solution is unique
 \begin{equation}\label{9.44}
 \bm h_k=-\|\bm g\|(T_k+\la_kI)^{-1}\bm e_1~~{\rm and}~~ \bm x_k=-\|\bm g\|Q_k(T_k+\la_kI)^{-1}\bm e_1.
 \end{equation}
Let $k_{\max}$ be the first number of iteration at which the symmetric Lanczos process breaks down, then we have $\bm x_{k_{\max}}=\bm x_{{opt}},~ \la_{k_{\max}}=\lambda_{opt}~~ {\rm and}~~ f(\bm x_{k_{\max}})=f(\bm x_{opt})$ \cite{10,5,4}.
We notice that the sequence $\{\la_k\}^{k_{\max}}_{k=0}$ of Lagrangian multipliers associated with $\eqref{9.25}$ is monotonically nondecreasing \cite{L}
\begin{equation}\label{eq9.42}
0\leq\la_1\leq \la_2\leq\cdots\leq\la_{k_{\max}}=\la_{opt}.
\end{equation}

\subsection{Some existing results}
Denote by
$$
A_{opt}= A+\la_{opt}I,~~\kappa=\frac{\alpha_1+\la_{opt}}{\alpha_n+\la_{opt}}~~{\rm and}~~\epsilon_k= \min_{\bm x\in \mathcal{K}_k(A,\bm g)}\|\bm x-\bm x_{opt}\|.
$$
Then we have from \cite[Theorem 4.10]{5} that
\begin{equation}\label{11.09}
\e_k=\|(I-Q_kQ^T_k)\bm x_{opt}\|\leq 2\D
\cdot {\Big(}\frac{\sqrt{\kappa}-1}{\sqrt{\kappa}+1}{\Big)}^{k}.
\end{equation}
Notice that this bound is slightly worse but is more concise than the one given in \cite[eq. (4.27)]{4}. Recently,
Z. Jia {\it et al}. \cite{5} established the following bound on $\|(A+\la_kI)\bm x_k+\bm g\|$:
\begin{theorem}\cite[Theorem 5.7]{5}
Suppose that $\|\bm x_{opt}\|=\|\bm x_k\|=\D$, $k=1,2,\ldots,k_{\max}$, and define $\bm r_k=(A+\la_{k}I)\bm x_k+\bm g$.
 Then asymptotically
\begin{equation}\label{10.27}
\|\bm r_k\|\leq
4\sqrt{\kappa}\|A_{opt}\|\D\cdot  {\Big(}\frac{\sqrt{\kappa}-1}{\sqrt{\kappa}+1}{\Big)}^{k}+{\Big(} \frac{4\eta_{k1}\D^2}{\|\bm g\|}+8\|A_{opt}\|\eta_{k2}\D^3
{\Big)}\cdot{\Big(}\frac{\sqrt{\kappa}-1}{\sqrt{\kappa}+1}{\Big)}^{2k}~~~~~~~
  \end{equation}
  with
  \begin{equation}\label{19.42}
  \begin{aligned}
\eta_{k1}= \frac{\|\bm g\|^2}{\D^2+\|\bm g\|^2\bm e_1^T(T_{k}+\la_{opt}I)^{-2}\bm e_1}\leq\frac{\|\bm g\|^2\|A_{opt}\|^2}{\|\bm g\|^2+\|A_{opt}\|^2\D^2},\\
\eta_{k2}= \frac{2}{\D^2+\|\bm g\|^2\bm e_1^T(T_{k}+\la_{opt}I)^{-2}\bm e_1}\leq\frac{2\|A_{opt}\|^2}{\|\bm g\|^2+\|A_{opt}\|^2\D^2}.
\end{aligned}
\end{equation}
\end{theorem}

In addition, Golud {\it et al}. \cite{G.S.} established another upper bound on $\|\bm r_k\|$:
\begin{theorem}\cite[Theorem 3.4]{G.S.}
The residual $\bm r_k=(A+\la_{k}I)\bm x_k+\bm g$ for the k-th iterate, $\bm x_k$, generated by the TRS \eqref{9.00} satisfies the
bound
\begin{equation}\label{10.28}
\|\bm r_k\| \leq\|\bm g\|{\Big(}\frac{2\beta_k\kappa_k}{\|T_k+\lambda_kI\|}{\Big)} {\Big(}\frac{\sqrt{\kappa_k}-1}{\sqrt{\kappa_k}+1}{\Big)}^{k-1},
\end{equation}
where $\kappa_k$ is the 2-condition number of $T_k\!+\!\la_kI$ and $\beta_k$ is the $(k,k\!+\!1)$-st entry of $T_{k+1}$.
\end{theorem}

Z. Jia {\it et al}. established the following {\it asymptotic} bound on the convergence of $\la_k$. That is, the bound holds only when $k$ is {\it sufficiently large}.
\begin{theorem}\cite[Theorem 4.11]{5}\label{14.144}
Suppose that $\|\bm x_k\|=\|\bm x_{opt}\|=\D$, $k=1,2,\ldots,k_{\max}$. Then asymptotically
\begin{equation}\label{11.15}
0\leq\la_{opt}-\la_k\leq {\Big (}\frac{4\eta_{k1}\D}{\|\bm g\|}+8\|A_{opt}\|\eta_{k2}\D^2{\Big )}{\Big(}\frac{\sqrt{\kappa}-1}{\sqrt{\kappa}+1}{\Big)}^{2k},
\end{equation}
where the constants $\eta_{k1}$ and $\eta_{k2}$ are defined by \eqref{19.42}.
\end{theorem}
\begin{remark}
{\it Theorem \ref{14.144} provides an ``asymptotic" bound, whose proof requires $k$ is sufficiently large \cite[pp.902]{5}. Therefore, it is interesting to establish ``non-asymptotic" bounds that hold in each step. In Section 4, we establish a ``non-asymptotic" bound which are closely related to the value of $s(\lambda_{opt})$ defined in \eqref{eq513}.}
\end{remark}

It was shown that \eqref{1} can be treated by solving the eigenvalue problem
on the $2n$-by-$2n$ matrix \cite{3,GGM,5}:
\begin{equation}\label{15.28}
M = \begin{pmatrix}
-A & \frac{\bm g\bm g^T}{\D^2}\\
I  &  -A
\end{pmatrix}\in \mathbb{R}^{2n\times 2n}.
\end{equation}
The following theorem establishes an important
relationship between the TRS solution $(\la_{opt},\bm x_{opt})$ and the rightmost eigenpair of $M$.
\begin{theorem}\label{16.0200}\cite{5}
Let $(\lambda_{opt}, \bm x_{opt})$ satisfy Theorem \ref{7.59} with $\|\bm x_{opt}\| = \D$.
Then the rightmost eigenvalue $\varpi_1$ of $M$ is real and simple, and $\la_{opt} = \varpi_1$.
  Let $\bm y = (\bm y^T_1,\bm y_2^T )^T$
be the corresponding unit length eigenvector of $M$ with $\bm y_1, \bm y_2 \in \mathbb{R}^{n}$, and suppose that $\bm g^T \bm y_2 \not = 0$. Then the unique TRS solution is
\begin{equation}\label{15.21}
\bm x_{opt} = -\frac{\D^2}{\bm g^T\bm y_2}\bm y_1.
\end{equation}
\end{theorem}
\begin{remark}
{\it It was shown that $\bm g^T\bm y_2=0$ only if TRS is in $\emph{hard case}$ \cite[Proposition 4.1]{3}. In other words, we have  $\bm g^T\bm y_2\neq 0$ in easy case.}
\end{remark}

In \cite{3}, Adachi \emph{et al}. proved that $\la_{opt}$ is the rightmost real eigenvalue of $M$. However, the proof of the uniqueness of $\la_{opt}$ was not given. For completeness, we give a proof here.

\begin{theorem}\label{16.0005}
Suppose that $\|\bm x_{opt}\|=\D$. Then
 trust-region subproblem \eqref{1} is in {easy case} if and only if $\la_{opt}$ is a simple eigenvalue of $M$.
\end{theorem}
\begin{proof}
On one hand, if TRS \eqref{1} is in \emph{easy case}, then $A_{opt}=A+\lambda_{opt}I\succ\bm O$, i.e., $\la_{opt}\in(-\alpha_n,\infty) $.
From \cite[eq. (14)--eq. (16)]{3}, we have
\begin{small}
\begin{equation}\label{20.355}
\det(M-\la I)=\det
\begin{pmatrix}
-I & A+\la I\\
A+\la I & -\frac{\bm g\bm g^T}{\D^2}
\end{pmatrix}
=(-1)^n\det(A+\la I)^2{\bigg (}1-\frac{\|(A+\la I)^{-1}\bm g\|^2}{\D^2}{\bigg )},
\end{equation}
where $\la\in(-\alpha_n,\infty)$.
\end{small}
Notice that $A_{opt}\succ\bm O$ and
$$
\|(A+\la I)^{-1}\bm g\|^2=\bm g^T(A+\la I)^{-2}\bm g=\sum^{n}_{i =1}\frac{(\bm u_i^T\bm g)^2}{(\alpha_i+\la)^2}.
$$
 Therefore,
$$
\frac{{\rm d} {\big ( }\det(M-\la I){\big ) }}{{\rm d} \la}{\Bigg \vert}_{\la=\la_{opt}}=(-1)^{n+2}\cdot\frac{2~ \det(A_{opt})^2}{\D^2}\cdot\left(
\sum^{n}_{i =1}\frac{(\bm u_i^T\bm g)^2}{(\alpha_i+\la_{opt})^3}\right)\neq 0.
$$
That is, $\la_{opt}$ is a simple root of $\det(M-\la I)=0$, i.e., $\la_{opt}$ is a simple eigenvalue of $M$.

On the other hand, if $\la_{opt}$ is a simple eigenvalue of $M$, we next prove that TRS \eqref{1} is in \emph{easy case}.
Recall that $\alpha_n$ is the smallest eigenvalue of $A$.
Suppose that the multiplicity of $\alpha_n$ is $s$, with $s>1$. Let
$$U_1=[\bm u_1,\bm u_2,\ldots,\bm u_{n-s}]~~{\rm and}~~U_2=[\bm u_{n-s+1},\bm u_{n-s},\ldots,\bm u_{n}].$$
If TRS \eqref{1} is in \emph{hard case}, i.e., $\la_{opt}=-\alpha_n$, then it follows that $ U^T_2\bm g= \bm 0$ and $\|(A-\alpha_n I)^{\dag}\bm g\|\leq\D$. Define
 $$\bm p(\la)=-\sum_{i=1}^{n-s+1}\frac{\bm u_i\bm u^T_i\bm g}{\alpha_i+\la},\quad\forall\la \in [-\alpha_n, \infty).$$

By \cite[eq. (14)--eq. (16)]{3},
\begin{equation}\label{15.510}
\det(M\!-\!\la I)\!=\! \frac{(-1)^{n}}{\D^2}\!\cdot\!\det
\begin{pmatrix}
\D^2\!-\!\|\bm p(\la)\|^2   &\! -\!\bm p(\la)^T &  \bm g^T\!+\!\bm p(\la)^T\!(A\!+\!\la I) \\
-\!\bm p(\la)  &\! -I \!& A\!+\!\la I\\
\bm g\!+\!(A\!+\!\la I)\bm p(\la) & A\!+\!\la I & \bm O
\end{pmatrix}.
\end{equation}
As $U_2^T\bm g=\bm 0$, we have $U_1U_1^T\bm g=\bm g$, and
\begin{align}\label{eq21180}
\bm g+(A+\la I)\bm p(\la)=\bm g-\sum_{i=1}^{n-s+1}\bm u_i\bm u^T_i\bm g=\bm g- U_1U_1^T\bm g=U_2U_2^T\bm g=\bm 0.
\end{align}
In addition, if $\la>-\alpha_n$, then $\bm p(\la)\in \mathcal{R}(A+\la I)=\mathbb{R}^n$. If $\la=-\alpha_n$, then
$$
\bm p(\la)=\bm p(-\alpha_n)\in\mathcal{R}(U_1)= \mathcal{R}(A-\alpha_n I)=\mathcal{R}(A+\la I).
$$

Thus, for any $\la\in [-\alpha_n,\infty)$, there is a vector $\bm q(\la)$ such that $\bm p(\la)=(A+\la I)\bm q(\la)$.
As a result,
\begin{equation}\label{eq2118}
\det
\begin{pmatrix}
\D^2\!-\!\|\bm p(\la)\|^2   & -\bm p(\la)^T &  \bm 0 \\
-\bm p(\la)  & -I & A\!+\!\la I\\
\bm 0 & A\!+\!\la I & \bm O
\end{pmatrix}
=(-1)^n \det(A\!+\!\la I)^2\cdot{\big (} \D^2\!-\!\|\bm p(\la)\|^2  {\big )}.
\end{equation}
It then follows from \eqref{15.510}--\eqref{eq2118} that
\begin{align*}
\det(M\!-\!\la I)\!= \!\det(A\!+\!\la I)^2{\Big [} 1\!-\!\frac{\|\bm p(\la)\|^2 }{\D^2}{\Big ]}
           \! =\!(\alpha_n\!+\!\la)^{2s}\left[\prod^{n-s+1}_{i=1}(\alpha_i\!+\!\la)^2{\Big (} 1\!-\!\frac{\|\bm p(\la)\|^2 }{\D^2}{\Big )}\right].
\end{align*}
Therefore, $\la_{opt}=-\alpha_n$ is a multiple root of $\det(M-\la I)=0$. In other words, $\la_{opt}$ is a multiple eigenvalue of $M$, a contradiction.
\end{proof}

Denote by
\begin{align}\label{19.54400}
\widetilde{Q}_k=\begin{pmatrix}
Q_k & \\
&  Q_k
\end{pmatrix}\in \mathbb{R}^{2n\times 2k}~~{\rm and}~~M_k=\widetilde{Q}_k^TM\widetilde{Q}_k=\begin{pmatrix}
-T_k  &  \frac{\|\bm g\|^2}{\D^2}\bm e_1\bm e^T_1\\
I_k    & -T_k
\end{pmatrix}.
\end{align}
Suppose that $\|\bm x_k\|=\D$. As $\la_k$ is the Lagrange multiplier of TRS \eqref{9.25}, it follows from Theorem \ref{16.0200} that $\la_k$ is the rightmost eigenvalue of $M_k$. Let $(\la_k,\bm z_k)$ be the rightmost eigenpair of $M_k$, where
$\bm z_k={\big(}(\bm z^{(k)}_{1})^T,(\bm z^{(k)}_{2})^T{\big)}^T$, with $\bm z^{(k)}_{1},\bm z^{(k)}_{2}\in \mathbb{R}^k$ and $\|\bm z_k\|=1$. Then
\begin{equation}\label{eq11.08}
\bm y_k=\widetilde{Q}_k\bm z_k={\big(}(Q_k\bm z^{(k)}_{1})^T,(Q_k\bm z^{(k)}_{2})^T{\big)}^T={\big (}(\bm y^{(k)}_{1})^T,(\bm y^{(k)}_{2})^T{\big)}^T
\end{equation}
is a Ritz vector of $M$ in the subspace $\widetilde{\mathcal{S}}_k=\mathcal{R}(\widetilde{Q}_k)$, which is an approximation to $\bm y$.

Let $Z_k$ be an orthnormal basis of the orthogonal complement of the subspace $span\{\bm z_k\}$, and let $C_k=Z^T_kM_kZ_k$. If $\la_{opt}$ is not an eigenvalue of $C_k$, we can define
 \begin{equation}\label{20.55}
sep(\lambda_{opt}, C_k)=\|(\la_{opt}I-C_k)^{-1}\|^{-1}.
\end{equation}
The following theorem depicts the distance between the approximation $\bm x_k$ and the optimal solution $\bm x_{opt}$:
\begin{theorem} \cite[Theorem 5.4]{5}\label{Thm4.2}
Suppose that $\|\bm x_{opt}\|\! = \!\| \bm x_k\| \!= \!\Delta$ and $sep(\lambda_{opt}, C_k)\!>\! 0$. Then
\begin{equation}\label{11.2300}
\sin\angle(\bm x_{opt},\bm x_k)\leq c_k\cdot\left(1\!+ \! \frac{\|M\|}{\sqrt{1-\sin^2\angle(\bm y, \widetilde{\mathcal{S}}_k)}\cdot sep(\la_{opt},C_k)}\right){\Big(}\frac{\sqrt{\kappa}-1}{\sqrt{\kappa}+1}{\Big)}^{k},
\end{equation}
where 
$$
c_k=2+\frac{16\|A_{opt}\|}{(\alpha_1-\alpha_n)^2(1-t^2)}{\Big(} 1+\frac{k+2}{\vert\ln t\vert}{\Big)} {\Big(}\frac{\sqrt{\kappa}-1}{\sqrt{\kappa}+1}{\Big)}^{2}~~{with}~~t=\frac{\sqrt{\kappa}-1}{\sqrt{\kappa}+1}.
$$
\end{theorem}
\begin{remark}\label{14.33}
From the differential mean value theorem, there is a constant $\nu\in [\sqrt{\kappa}-1,\sqrt{\kappa}+1]$, such that
\begin{equation}\label{eq45}
(1-t)|\ln t|
=\frac{2}{\sqrt{\kappa}+1}{\big(}\ln(\sqrt{\kappa}+1)-\ln(\sqrt{\kappa}-1)
{\big)}=\frac{4}{(\sqrt{\kappa}+1)\nu}\leq\frac{4}{\kappa-1},
\end{equation}
and thus $\frac{1}{(1-t)\vert\ln t\vert}=\mathcal{O}(\kappa)$. Moreover, we have that
\begin{align}\label{17.06}
\frac{1}{(\alpha_1\!-\!\alpha_n)}\cdot\frac{\sqrt{\kappa}\!-
\!1}{\sqrt{\kappa}\!+\!1}=\frac{\kappa}{\|A_{opt}\|(\sqrt{\kappa}\!+\!1)^2}
\!=\!\mathcal{O}\Big(\frac{1}{\|A_{opt}\|}\Big).
\end{align}
Therefore, it follows from \eqref{eq45} and \eqref{17.06}
that $c_k=\mathcal{O}{\big(}\frac{k\cdot\kappa}{\|A_{opt}\|}{\big)}$.
In summary,
\begin{equation}\label{eq47}
\widetilde{c}_k= c_k\cdot\left(1+ \frac{\|M\|}{\sqrt{1-\sin^2\angle(\bm y, \widetilde{\mathcal{S}}_k)}\cdot sep(\la_{opt},C_k)}\right)=\mathcal{O}\left(\frac{ k \kappa}{sep(\la_{opt},C_k)}\right),
\end{equation}
i.e.,
$$
\sin\angle(\bm x_{opt},\bm x_k)\leq \mathcal{O}\left(\frac{ k \kappa}{sep(\la_{opt},C_k)}\right)
{\Big(}\frac{\sqrt{\kappa}-1}{\sqrt{\kappa}+1}{\Big)}^{k}.
$$
In other words, \eqref{11.2300} shows that the upper bound of $\sin\angle(\bm x_{opt},\bm x_k)$ is closely related to $sep(\la_{opt},C_k)$, $\kappa$ and $k$.
However, there is no guarantee that $sep(\la_{opt},C_k)$ is uniformly lower bounded, both in theory and in practice.
Consequently, the righthand side of \eqref{11.2300} can be a poor estimation to the distance between $\bm x_k$ and $\bm x_{opt}$.
\end{remark}

\section{Improved upper bounds on $\bm{\|(A+ \la_{k} I) x_k+g\|}$ and ${\bf sin\angle}\bm{(x_k, x_{opt})}$}
 In this section, we give some new upper bounds on ${\|(A+ \la_{k} I) {\bm x_k}+{\bm g}\|}$ and ${\sin\angle}(\bm x_k, \bm x_{opt})$.
Let $\mathcal{P}_{i}$ be the set of polynomials whose degrees are no higher than the positive integer $i$. Let $\mathcal{C}_i(x)\in \mathcal{P}_i$ be the Chebyshev polynomial define as
\begin{equation}\label{def2004}
\mathcal{C}_{i}(x)=\left\{
\begin{array}{ll}
\cos(i\arccos x)                     &  {\rm for}~~ \vert x\vert<1,\\
\specialrule{0em}{0.8ex}{0.8ex}
\frac{1}{2}{\big(}(x+\sqrt{x^2-1})^i+(x+\sqrt{x^2-1})^{-i} {\big)}~~ & {\rm for}~~\vert x\vert\geq 1.
\end{array}\right.
\end{equation}
To establish our bounds, we need the following result.
\begin{lemma}
Given a matrix $H\!\succ\!{\bm O}$. Let $\kappa_H\!=\!\frac{\la_{\max}(H)}{\la_{\min}(H)}$ be the 2-condition number of $H$.
\begin{itemize}
\item[\rm (i)]
\cite[Lemma 1,~Lemma 4]{110}\label{9.210}
  Suppose that $\varphi_i\in \mathcal{P}_i$  satisfies
\begin{equation}\label{1.200}
1\!-\!x\varphi_i(x)\!=\!2{\big(} r^{i+1}\!+\!r^{-(i+1)}{\big)}\!^{-1}\cdot\mathcal{C}_{i+1}
{\bigg(}\frac{\kappa_H\!+\!1}{\kappa_H\!-\!1}\!-\!\frac{2x}
{\la_{\max}(H)\!-\!\la_{\min}(H)} {\bigg)},
\end{equation}
where $r=\frac{\sqrt{\kappa_H}-1}{\sqrt{\kappa_H}+1}$.
Then, we have that
$\|I-H\varphi_i(H)\|\leq 2r^{i+1}.$
\item[\rm (ii)] Suppose that $\widetilde{H}\!$ is symmetric and $[\la_{\min}(\widetilde{H}),\la_{\max}(\widetilde{H})]\!\subseteq\! [\la_{\min}({H}),\la_{\max}({H})]$, then $\|I\!-\!\widetilde{H}\varphi_i(\widetilde{H})\|\!\leq\! 2r^{i+1}$. 
\end{itemize}
\end{lemma}
\begin{proof}
We only need to prove (ii).
Let $\widetilde{H} = P\Gamma P^T$ be the eigendecomposition,
where $P$ is orthonormal, and
$
\Gamma=diag(\gamma_1,\gamma_2,\ldots,\gamma_m)
$
with $\la_{\max}(\widetilde{H})=\gamma_1\geq\gamma_2
\geq\cdots\geq\gamma_m=\la_{\min}(\widetilde{H})$ being the eigenvalues of $\widetilde{H}$. Thus,
\begin{align*}
\|I-\widetilde{H}\varphi_i(\widetilde{H})\|
&=\|I-\Gamma\varphi_i(\Gamma)\|=\max_{\ell=1,\ldots,m}
\vert1-\gamma_\ell\varphi_i(\gamma_\ell)\vert\\
&=2{\big(} r^{i+1}\!+\!r^{-(i+1)}{\big)}\!^{-1}
\cdot\underset{=\varrho}{\underbrace{\max_{\ell=1,\dots,m}
\left\vert\mathcal{C}_{i+1}
{\bigg(}\frac{\kappa_H\!+\!1}{\kappa_H\!-\!1}\!-\!\frac{2\gamma_\ell}
{\la_{\max}(H)\!-\!\la_{\min}(H)} {\bigg)}\right\vert}}.
\end{align*}
{As $\gamma_\ell\in[\la_{\min}(H),\la_{\max}(H)]$, we have that
$$
\left\vert\frac{\kappa_H+1}{\kappa_H-1}-\frac{2\gamma_\ell}
{\la_{\max}(H)-\la_{\min}(H)}\right\vert\leq 1,
$$
and it follows from \eqref{def2004} that $\varrho\leq 1$.}
Hence, {\small$\|I\!-\!\widetilde{H}\varphi_i(\widetilde{H})\|\leq2{\big(} r^{i+1}\!+\!r^{-(i+1)}{\big)}\!^{-1}\leq 2r^{i+1}$}, which completes the proof.
\end{proof}

We are ready to establish some improved bounds on $\|(A+\la_kI)\bm x_k+\bm g\|$ and $\sin\angle(\bm x_{opt},\bm x_k)$:
\begin{theorem}\label{9.051}
Suppose that $\|\bm x_{opt}\|=\|\bm x_k\|=\D$. Then
\begin{align}
&\|(A+\la_kI)\bm x_k+\bm g\|\leq \min\{\xi_1,\xi_2\},\label{10.571}\\
&\sin\angle(\bm x_{opt},\bm x_k)\leq2\sqrt{\kappa{\Big(}1+\frac{{\epsilon_k^2}}{\D^2}{\Big)}}
\cdot{\Big(}\frac{\sqrt{\kappa}\!-\!1}{\sqrt{\kappa}\!+\!1}{\Big)}^{k},
\label{16.0600}
\end{align}

where
$$
\xi_1=2\|A_{opt}\|\sqrt{\D^2+\e_k^2}\cdot{\Big(}\frac{\sqrt{\kappa}-1}{\sqrt{\kappa}+1}{\Big)}^{k}
~~{and}~~
\xi_2=2\beta_k\D\cdot{\Big(}\frac{\sqrt{{\kappa}_k}-1}{\sqrt{{\kappa}_k}+1}{\Big)}^{k-1}.
$$
\end{theorem}
\begin{proof}First,
for any $\|\bm x\| = \D$, we have that
\begin{align}\label{2.39}
0\leq f(\bm x)-f(\bm x_{opt})&=\frac{1}{2}\left(\bm x^TA\bm x-\bm x_{opt}^TA\bm x_{opt}\right)
+\bm g^T(\bm x-\bm x_{opt})\nonumber\\
         &=\frac{1}{2}\left(\bm x^TA_{opt}\bm x-\bm x_{opt}^TA_{opt}\bm x_{opt}\right)
         -\bm x_{opt}^TA_{opt}(\bm x-\bm x_{opt})\nonumber\\
         &=\frac{1}{2}(\bm x_{opt}-\bm x)^TA_{opt}(\bm x_{opt}-\bm x).
\end{align}
Let ${\bm y}_k=\D\cdot\frac{Q_kQ^T_k\bm x_{opt}}{\|Q_kQ^T_k\bm x_{opt}\|}\in \mathcal{K}_k(A,\bm g)$, then
it follows from \eqref{9.000} and \eqref{2.39} that
\begin{align}
f({\bm x}_k)\!-\!f(\bm x_{opt})\!\leq\! f{(}{\bm y}_k{)}\!-\!f(\bm x_{opt})
\leq\frac{1}{2}\|A_{opt}\|\|\bm
x_{opt}\!-\!\bm y_k\|^2.\label{10.24}
\end{align}
Moreover,
\begin{align}
\|\bm x_{opt}\!-\!\bm y_k\|^2 &=\! \|(I-Q_kQ_k^T)\bm x_{opt}+Q_kQ_k^T\bm x_{opt}-\bm y_k\|^2\nonumber\\
                          &=\!\|(I-Q_kQ_k^T)\bm x_{opt}\|^2+\|Q_kQ_k^T\bm x_{opt}-\bm y_k\|^2\nonumber\\
                          &=\! \|(I-Q_kQ_k^T)\bm x_{opt}\|^2+\|Q_kQ_k^T\bm x_{opt}\|^2+\|\bm y_k\|^2-2\bm y^T_kQ_kQ_k^T\bm x_{opt}\nonumber\\
                          &=\!\|(I\!-\!Q_kQ_k^T)\bm x_{opt}\|^2\!+\!{\big(}\D\!-\!\|Q_kQ_k^T\bm x_{opt}\|{\big)}^2\!=\!\e_k^2\!+\!
                          {\big(}\D\!-\!\|Q_kQ_k^T\!\bm x_{opt}\|{\big)}^2.\nonumber
\end{align}
Now we consider $\D-\|Q_kQ_k^T\bm x_{opt}\|$. As $\|\bm x_{opt}\|=\D$ and $\|Q_kQ_k^T\bm x_{opt}\|\!\leq\!\|\bm x_{opt}\|\!=\!\D$, we have that
\begin{align*}
\D\!-\!\|Q_kQ_k^T\!\bm x_{opt}\|\!=\!\frac{\D^2\!-\!\D\|Q_k^T\!\bm x_{opt}\|}{\D}\!\leq\!\frac{\D^2\!-\!\|Q_k^T\!\bm x_{opt}\|^2}{\D}
\!=\!\frac{\|(I\!-\!Q_kQ_k^T)\bm x_{opt}\|^2}{\D}\!=\!\frac{\e_k^2}{\D}.
\end{align*}
As a result,
\begin{equation}\label{eq10.29}
0\!\leq\!f(\bm x_k)\!-\!f(\bm x_{opt})\!\leq\!\frac{1}{2}\|A_{opt}\|{\Big (}1\!+\!\frac{\epsilon_k^2}{\D^2} {\Big )}
\!\cdot\! \epsilon_k^2\overset{\eqref{11.09}}{\leq}
2\|A_{opt}\|{\Big(}\D^2\!+\!{\epsilon_k^2}{\Big)}
\cdot{\Big(}\frac{\sqrt{\kappa}\!-\!1}{\sqrt{\kappa}\!+\!1}{\Big)}^{2k}.
\end{equation}

Notice that
\begin{align}
A_{opt}(\bm x_k-\bm x_{opt})=(A+\la_{opt}I)\bm x_k+\bm g
                                   =(A+\la_{k}I)\bm x_k+\bm g+(\la_{opt}-\la_{k})\bm x_k.\label{11.59}
\end{align}
By \eqref{11.10},
$$\|A_{opt}(\bm x_k-\bm x_{opt})\|^2=\|(A+\la_{k}I)\bm x_k+\bm g\|^2+(\la_{opt}-\la_{k})^2\|\bm x_k\|^2,$$
then
\begin{align}\label{eq2037}
&\|(A+\la_{k}I)\bm x_k+\bm g\|^2\leq\|A_{opt}(\bm x_k-\bm x_{opt})\|^2\leq\|A^{\frac{1}{2}}_{opt}\|^2\cdot\|A^{\frac{1}{2}}_{opt}(\bm x_k-\bm x_{opt})\|^2\nonumber\\
 =&\|A_{opt}\|\cdot(\bm x_k\!-\!\bm x_{opt})^T\!A_{opt}(\bm x_k\!-\!\bm x_{opt})\!\overset{\eqref{2.39}}{=}\!2\|A_{opt}\|\cdot\big(f(\bm x_k)\!-\!f(\bm x_{opt})\big)\!\overset{\eqref{eq10.29}}{\leq}\! \xi^2_1.
\end{align}

Second, for any $\varphi_{k-2}\in \mathcal{P}_{k-2}$, we have from \eqref{7.55} that
\begin{align}\label{eq2042}
\|(A\!+\!\la_{k}I)\bm x_k\!+\!\bm g\|\! = \!\beta_k\vert\bm e_k^T\!\bm h_k\vert&\!\leq\!\beta_k {\Big(}{\big\vert}\bm e_k^T{\big(}I-(T_k+\la_{k}I)\varphi_{k-2}(T_k+\la_{k}I){\big)}\bm h_k{\big\vert}\nonumber\\
&~~~~~~~~~~+ {\big\vert}\bm e_k^T(T_k+\la_{k}I)\varphi_{k-2}(T_k+\la_{k}I)\bm h_k{\big\vert}{\Big)}.
\end{align}
Recall that $\varphi_{k-2}\in\mathcal{P}_{k-2}$ and $T_k+\la_{k}I$ is a symmetric and tridiagonal matrix.
Thus, $\bm e_k^T\varphi_{k-2}(T_k+\la_{k}I)\bm e_1=0$ for any $ \varphi_{k-2}\in \mathcal{P}_{k-2}$. Consequently,
\begin{align*}
\bm e_k^T\!(T_k\!+\!\la_{k}I)\cdot\varphi_{k-2}(T_k\!+\!\la_{k}I)\bm h_k&=-\|\bm g\|\bm e_k^T\!(T_k\!+
\!\la_{k}I)\cdot\varphi_{k-2}(T_k\!+
\!\la_{k}I)(T_k\!+\!\la_{k}I)^{-1} \bm e_1\\
&=-\|\bm g\|\bm e_k^T\varphi_{k-2}(T_k+\la_{k}I)\bm e_1=0.
\end{align*}
Specially, if we choose $\varphi_{k-2}$ that satisfies \eqref{1.200}, it then follows from Lemma \ref{9.210} that
\begin{equation}\label{12.15}
\|(A+\la_{k}I)\bm x_k+\bm g\|\overset{\eqref{eq2042}}{\leq} \beta_k\|\bm h_k\|\|I-(T_k+\la_{k}I)\varphi_{k-2}(T_k+\la_{k}I)\|\leq \xi_2.
\end{equation}
A combination of the above equation with \eqref{eq2037} gives \eqref{10.571}.

Finally, we notice that
\begin{align*}
\left(\sin\angle(\bm x_{opt},\bm x_k)\right)^2&=\min_{\chi\in \mathbb{C}}{\Big \|}\frac{\bm x_{opt}}{\D}- \chi\cdot \bm x_k{\Big \|}^2\leq
{\Big \|}\frac{\bm x_{opt}}{\D}- \frac{\bm x_k}{\D}{\Big \|}^2=\frac{\|\bm x_{opt}-\bm x_k\|^2}{\D^2}\\
        &\leq \frac{(\bm x_{opt}-\bm x_k)^TA_{opt}(\bm x_{opt}-\bm x_k)}{\D^2(\alpha_n+\la_{opt})}\overset{\eqref{2.39}}{=}\frac{2 \big(f(\bm x_{k})-f(\bm x_{opt})\big)}{\D^2(\alpha_n+\la_{opt})}\\
        &\overset{\eqref{eq10.29}}{\leq}4\kappa{\Big(}1+\frac{{\epsilon_k^2}}{\D^2}{\Big)}
\cdot{\Big(}\frac{\sqrt{\kappa}\!-\!1}{\sqrt{\kappa}\!+\!1}{\Big)}^{2k},
\end{align*}
which completes the proof.
\end{proof}
\begin{remark}
{\it First, the upper bound established in \eqref{10.27} is about $\sqrt{\kappa}$ times larger than $\xi_1$,
and our new upper bound \eqref{10.571} improves the one given
in \eqref{10.27}, especially when the condition number $\kappa$ is large.
Moreover, we notice that
$$
\D = \|\bm x_k\|=\|\bm h_{k}\|=\|\bm g\|\|(T_k+\la_kI)^{-1}\bm e_1\|\leq
{\|(T_k+\la_kI)^{-1}\|\|\bm g\|}=\frac{\|\bm g\|\kappa_k}{\|T_k+\lambda_kI\|},
$$
and our bound \eqref{10.571} is no worse than \eqref{10.28}.
Indeed, compared with the bounds \eqref{10.27} and \eqref{10.28}, in our new upper bound  \eqref{10.571}, there are no
condition numbers such as $\kappa$ or $\kappa_k$ in the coefficients before ${\big(}\frac{\sqrt{{\kappa}_k}-1}{\sqrt{{\kappa}_k}+1}{\big)}^{k-1}$ and ${\big(}\frac{\sqrt{{\kappa}}-1}{\sqrt{{\kappa}}+1}{\big)}^{k}$, respectively.
Hence, we improved the upper bounds of Z. Jia {\it et al}. and Gould et al. substantially.

Second,
our new upper bound \eqref{16.0600} improved \eqref{11.2300} substantially.
More precisely, compared with  \eqref{11.2300}, \eqref{16.0600} indicates that the factor before $\big(\frac{\sqrt{\kappa}-1}{\sqrt{\kappa}+1}{\big)}^{k}$ is in the order of $\sqrt{\kappa}$, which is free of $sep(\la_{opt},C_k)$ and $k$. In other words, Theorem \ref{Thm4.2} shows that the convergence of $\bm x_k$ is influenced by the factor $sep(\la_{opt},C_k)$ \cite{5}, while our bound \eqref{16.0600} reveals that it has nothing to do with $sep(\la_{opt},C_k)$ at all.}
\end{remark}

\section{Non-asymptotic error estimates for $\bm{\la_{{opt}}-\la_{{k}}}$}
In this section, we focus on the convergence of $\la_k$.
To do this, we first need the following three lemmas.
\begin{lemma}\label{lem19.59}
Let $t=\frac{\sqrt{\kappa}-1}{\sqrt{\kappa}+1}$, and denote by
$$
g_{i}(x)=2\big(t^{i}
 + t^{-i} \big)^{-1}\cdot\mathcal{C}_{i}{\bigg(}\frac{\kappa+1}{\kappa-1}-\frac{2x}{\alpha_1-\alpha_n} {\bigg)},
$$
where $\mathcal{C}_{i}(\cdot)$ is defined in \eqref{def2004}. Then we have that
$$
0\leq -\frac{{\rm d}g_{i}}{{\rm d}x} \bigg\vert_{x=0}\leq \frac{i\sqrt{\kappa}}{\|A_{opt}\|}~~{ and}~~\frac{{\rm d}^2g_{i}}{{\rm d}x^2} \bigg\vert_{x=0}\geq 0,\quad i=1,2,\ldots
$$
\end{lemma}
\begin{proof}
Let
$\ell=\frac{\kappa+1}{\kappa-1}$, then
\begin{equation}\label{20.37000}
\frac{{\rm d}g_i}{{\rm d}x} \bigg\vert_{x=0}=-\frac{4 ~\big(t^{i}
 + t^{-i} \big)^{-1}}{\alpha_1-\alpha_n}\cdot\frac{{\rm d}\mathcal{C}_{i}}{{\rm d}x}{\bigg \vert}_{x=\ell}.
\end{equation}
Since
$$
\left\vert{\Big(}x-\sqrt{x^2-1}{\Big)}^i {\Big(} 1-\frac{x}{\sqrt{x^2-1}} {\Big )}
\right\vert\leq{\Big(}x+\sqrt{x^2-1}{\Big)}^i {\Big(} 1+\frac{x}{\sqrt{x^2-1}} {\Big)},\quad \forall x>1,
$$
we have that
\begin{small}
\begin{align}\label{20.2700}
0\leq\frac{{\rm d}\mathcal{C}_{i}}{{\rm d}x}&={\frac{i}{2}}{\Bigg(} {\Big(}x+\sqrt{x^2-1}{\Big)}^{i-1} {\Big(} 1+\frac{x}{\sqrt{x^2-1}}
{\Big)}+{\Big(}x-\sqrt{x^2-1}{\Big)}^{i-1} {\Big(} 1-\frac{x}{\sqrt{x^2-1}} {\Big)} {\Bigg)}\nonumber\\
&\leq
{\frac{i}{2}}{\Big(}x+\sqrt{x^2-1}{\Big)}^{i-1} {\Big(} 1+\frac{x}{\sqrt{x^2-1}}
{\Big)},\quad \forall x>1.
\end{align}
\end{small}
Moreover,
\begin{equation*}
\ell+\sqrt{\ell^2-1}=\frac{\sqrt{\kappa}+1}{\sqrt{\kappa}-1},
~~1+\frac{\ell}{\sqrt{\ell^2-1}}=\frac{(\sqrt{\kappa}+1)^2}{2\sqrt{\kappa}}
~~{\rm and}~~\big(t^{i}
 + t^{-i} \big)^{-1}\leq t^{i}.
\end{equation*}
So it follows from \eqref{20.2700} that
\begin{align*}
0\leq\big(t^{i}
 + t^{-i} \big)^{-1}\cdot\frac{{\rm d}\mathcal{C}_{i}}{{\rm d}x}{\bigg \vert}_{x=l_1}\leq
{\frac{i}{4}}\cdot\frac{(\sqrt{\kappa}+1)^2}{\sqrt{\kappa}}\cdot {\Big (}\frac{\sqrt{\kappa}-1}{\sqrt{\kappa}+1}{\Big)}.
\end{align*}
By \eqref{20.37000},
\begin{align}
0\leq-\frac{{\rm d}g_{i}}{{\rm d}x} \bigg\vert_{x=0}
\leq\frac{i}{\alpha_1-\alpha_n}\cdot
\frac{(\sqrt{\kappa}+1
)^2}{\sqrt{\kappa}}\cdot{\Big(}\frac{\sqrt{
\kappa}-1}{\sqrt{\kappa}+1}{\Big)}=
\frac{i\sqrt{\kappa}}{\|A_{opt}\|},\label{8.27}
\end{align}
where the last equality is from \eqref{17.06}.

In addition, it follows from \eqref{20.37000} that
$$
\frac{{\rm d}^2{g}_{i}}{{\rm d}x^2}{\bigg \vert}_{x=0}= \frac{8 ~\big(t^{i}
 + t^{-i} \big)^{-1}}{(\alpha_1-\alpha_n)^2}\cdot\frac{{\rm
d}^2\mathcal{C}_{i}}{{\rm d}x^2}{\bigg\vert}_{x=\ell}.
$$
From \cite[eq. (13)]{Jia}, we obtain
$$
\frac{{\rm
d}^2\mathcal{C}_{i}}{{\rm d}x^2}{\bigg\vert}_{x=\ell}=\sum_{s=0}^{i-2}\vartheta_{s}
~\mathcal{C}_{s}\left(\ell\right)~~{\rm with}~~\vartheta_{s}\geq 0,~~s=0,1,\ldots,i-2.
$$
As $\ell>1$, we have $\mathcal{C}_i(\ell)\geq 0$, and $\frac{{\rm d}^2{g}_{i}}{{\rm d}x^2}{\big\vert}_{x=0}\geq 0$.
\end{proof}
\begin{lemma}\label{lem14.43}
For $j=1,2,\ldots,2k-2$, we have that
$$
\bm g^TA_{opt}^j\bm g=\|\bm g\|^2\cdot\bm e^T_1(T_k+\la_{opt}I)^j\bm e_1.
$$
\end{lemma}
\begin{proof}
We denote by $T_{k,opt}= T_k+\la_{opt}I$ for the sake of simplicity.
For $j=1,2,\ldots,k-1$, we have from \cite[eq. (A.9)]{G.S.} that
$$
A_{opt}^j Q_k=Q_kT_{k,opt}^j +\beta_k \cdot \sum^{j-1}_{i=0}A_{opt}^{j-i-1}\bm q_{k+1}\cdot\bm e_k^TT_{k,opt}^i.
$$
As $T_{k,opt}$ is a symmetric and tridiagonal matrix,  $\bm e_k^TT_{k,opt}^i\bm e_1=0$ for $0\leq i\leq k-2$. Recall that $Q_k\bm e_1=\frac{\bm g}{\|\bm g\|}$, so we obtain
\begin{eqnarray}\label{eq2129}
\left\{
\begin{array}{ll}
Q_kT_{k,opt}^j\bm e_1=A_{opt}^j Q_k\bm e_1=\frac{A_{opt}^j\bm g }{\|\bm g\|},\\
\specialrule{0em}{0.8ex}{0.8ex}
\bm e_1^TT_{k,opt}^j\bm e_1=\bm e_1^TQ^T_k A_{opt}^j Q_k\bm e_1=\frac{\bm g^T A_{opt}^j\bm g}{\|\bm g\|^2},
\end{array}\right.~~~~j=1,2,\ldots,k-1.
\end{eqnarray}
For $j=k,k+1,\ldots,2k-2$, we have $j-k+1=1,2,\ldots,k-1$, and
\begin{align*}
\frac{\bm g^TA_{opt}^j\bm g}{\|\bm g\|^2}=&\bm e_1^TQ^T_k A_{opt}^j Q_k\bm e_1=(A_{opt}^{j-k+1}Q_k\bm e_1)^T A_{opt}^{k-1} Q_k\bm e_1\\
=&\left(Q_kT_{k,opt}^{j-k+1}\bm e_1\right)^T   Q_kT_{k,opt}^{k-1}\bm e_1
= \bm e_1^TT_{k,opt}^{j}\bm e_1,
\end{align*}
where the third equality follows from the first equation of \eqref{eq2129}.
\end{proof}

Let
\begin{equation}\label{19.56}
\widetilde{\bm h}_k=-\|\bm g\|(T_k+\la_{opt} I)^{-1}\bm e_1~~~{\rm and}~~~\widetilde{\bm x}_k= Q_k\widetilde{\bm h}_k.
\end{equation}
Notice that
\begin{equation}\label{eq1637}
A_{opt}\widetilde{\bm x}_k+\bm g=A_{opt}Q_k\widetilde{\bm h}_k+\|\bm g\|Q_k\bm e_1\overset{\eqref{21.18}}{=} (\beta_k\bm e_1^T\widetilde{\bm h}_k)\cdot\bm q_{k+1}~\bot~\mathcal{K}_k(A,\bm g),
\end{equation}
then we have that
\begin{lemma}\label{22.22}
Let $\widetilde{\bm x}_{k}$ be defined in \eqref{19.56}, $\|\bm x_{opt}\|_{A_{opt}}=\sqrt{\bm x_{opt}^T{A_{opt}}\bm x_{opt}}$  and $k=1,2,\ldots,k_{\max}$. Then
\begin{eqnarray}\label{eq8.26}
0\leq \frac{\|\bm x_{opt}\|^2-\|\widetilde{\bm x}_k\|^2}{\D^2}\leq 4\left(1+\frac{2k\sqrt{\kappa}}{\|A_{opt}\|}\left(\frac{\|\bm x_{opt}\|_{A_{opt}}}{\D}\right)^2\right)
{\Big(}\frac{\sqrt{\kappa}-1}{\sqrt{\kappa}+1}{\Big)}^{2k}.~~~~~~~~
\end{eqnarray}
\end{lemma}
\begin{proof}
From \eqref{8.18} and the facts that $\la_{opt}\geq\la_k$ and $T_k+\la_{opt}I=Q_k^TA_{opt}Q_k\succ\bm O$, we obtain 
\begin{equation}\label{17.0003}
0\leq\|\bm x_{opt}\|^2-\|\widetilde{\bm x}_k\|^2=\bm x_{opt}^T\bm x_{opt}-\widetilde{\bm x}_{k}^T\widetilde{\bm x}_{k}=\bm g^TA_{opt}^{-2}\bm g-\|\bm g\|^2\cdot\bm e_1^TT_{k,opt}^{-2}\bm e_1,
\end{equation}
where $T_{k,opt}= T_k+\la_{opt}I$.
Suppose that $\varphi_{2k-1}(x)\in \mathcal{P}_{2k-1}$ satisfies \eqref{1.200}, and let
$$
p_{2k-2}(x)=\frac{\varphi_{2k-1}(x)-\widetilde{a}_0}{x}\in\mathcal{P}_{2k-2},~~{\rm with}~~x\in[\alpha_n+\la_{opt},\alpha_1+\la_{opt}],
$$
where $\widetilde{a}_0$ is the constant term of $ \varphi_{2k-1}(x)$. According to Lemma \ref{lem14.43},
$$
\bm g^Tp_{2k-2}(A_{opt})\bm g=\|\bm g\|^2\cdot\bm e_1^Tp_{2k-2}(T_{k,opt})\bm e_1,
$$
so we obtain
\begin{small}
\begin{align}
    \frac{\bm g^TA_{opt}^{-2}\bm g}{\|\bm g\|^2}-\bm e^T_1T_{k,opt}^{-2}\bm e_1=\frac{\bm g^T\big[A_{opt}^{-2}-p_{2k-2}(A_{opt})\big]\bm g}{\|\bm g\|^2}-\bm e^T_1\big[T_{k,opt}^{-2}-p_{2k-2}(T_{k,opt})\big]\bm e_1.\label{eq1522}
\end{align}
\end{small}

Now we consider $\bm g^T\big[A_{opt}^{-2}-q_{2k-2}(A_{opt})\big]\bm g$. By the definition of $p_{2k-2}(x)$, we get
\begin{align}
    &\bm g^T\big[A_{opt}^{-2}-p_{2k-2}(A_{opt})\big]\bm g
=\bm g^TA^{-1}_{opt}\big[I-A_{opt}^2p_{2k-2}(A_{opt})\big]A^{-1}_{opt}\bm g\nonumber\\
=&\bm x^T_{opt}\big[I-A_{opt}^2p_{2k-2}(A_{opt})-\widetilde{a}_0A_{opt}
\big]\bm x_{opt}+\widetilde{a}_0\bm g^TA^{-1}_{opt}\bm g\nonumber\\
=&\bm x^T_{opt}\big[I-A_{opt}\varphi_{2k-1}(A_{opt})\big] \bm x_{opt}-\widetilde{a}_0\bm g^T\bm x_{opt}.\label{eq1527}
\end{align}
Similarly, we can prove that
\begin{align}\label{eq1535}
\bm e_1^T\big[T_{k,opt}^{-2}-p_{2k-2}(T_{k,opt})\big]\bm e_1
=\widetilde{\bm h}_k^T\big[I-T_{k,opt}\varphi_{2k-1}(T_{k,opt})\big]
\widetilde{\bm h}_k-\frac{\widetilde{a}_0\cdot\bm e_1^T\widetilde{\bm h}_{k}}{\|\bm g\|},
\end{align}
where ${\widetilde{\bm h}}_k$ is defined in \eqref{19.56}.

From the Poincar\'{e} separation theorem \cite[Corollary 4.3.37]{H.B}, we have
$\alpha_n+\la_{opt}\leq\la_{\min}(T_k)+\la_{opt}
\leq\la_{\max}(T_k)+\la_{opt}\leq \alpha_1+\la_{opt}.$
Hence,
\begin{align}\label{eq15,35}
\|I-T_{k,opt}\varphi_{2k-1}(T_{k,opt})\|\overset{\rm Lem. \ref{9.210} (ii)}{\leq}
2\left(\frac{\sqrt{\kappa}-1}{\sqrt{\kappa}+1}\right)^{2k}.
\end{align}
It follows from \eqref{eq1522}--\eqref{eq15,35} and Lemma \ref{9.210} (i) that
\begin{align}
    &\bm g^TA_{opt}^{-2}\bm g-\|\bm g\|^2\cdot\bm e^T_1(T_{k}+\la_{opt}I)^{-2}\bm e_1\nonumber\\
  =&\bm x^T_{opt}\big[I-A_{opt}\varphi_{2k-1}(A_{opt})\big]
  \bm x_{opt}-\|\bm g\|^2\cdot\widetilde{\bm h}^T_k\big[I-T_{k,opt}\varphi_{2k-1}(T_{k,opt})\big]
\widetilde{\bm h}_k\nonumber\\
   &~~+\widetilde{a}_0(\|\bm g\|\cdot\bm e_1^T\widetilde{\bm h}_{k}-\bm g^T\bm x_{opt})\nonumber\\
   \leq&4\D^2{\Big(}\frac{\sqrt{\kappa}-1}{\sqrt{\kappa}+1}{\Big)}^{2k}+\widetilde{a}_0(\|\bm g\|\cdot\bm e_1^T\widetilde{\bm h}_{k}-\bm g^T\bm x_{opt}),\label{eq1529}
\end{align}
where the last inequality is from \eqref{eq15,35} and $\|\widetilde{\bm h}_k\|=\|\widetilde{\bm x}_k\|\leq\D$; refer to \eqref{17.0003}.

Moreover, by \cite[Theorem 3.1.1]{8},
\begin{subequations}
\begin{align}
\|\bm g\|\!\cdot\!\bm e_1^T\widetilde{\bm h}_{k}\!&-\!\bm g^T\!\bm x_{opt}\!
=\!\|\bm g\|\!\cdot\!\bm e_1^T\!Q^T_k\!Q_k\widetilde{\bm h}_{k}\!-\!\bm g^T\!\bm x_{opt}\!
=\!\bm g^T\!\widetilde{\bm x}_k\!-\!\bm g^T\!\bm x_{opt}\!
=\!-\bm x_{opt}^T\!(A_{opt}\widetilde{\bm x}_{k}\!+\!\bm g)\nonumber\\
&\overset{\eqref{eq1637}}{=}\!(\widetilde{\bm x}_{k}\!-\!\bm x_{opt})^T\!(A_{opt}\widetilde{\bm x}_{k}\!+\!\bm g)\! =\!(\widetilde{\bm x}_{k}\!-\!\bm x_{opt})^T\!A_{opt}(\widetilde{\bm x}_{k}\!-\!\bm x_{opt})\label{1932001}\\
&\leq4~\|\bm x_{opt}\|^2_{A_{opt}}\cdot \Big( \frac{\sqrt{\kappa}-1}{\sqrt{\kappa}+1}\Big)^{2k}.\label{1932000}
\end{align}
\end{subequations}
Recall that $\widetilde{a}_0$ is the constant term of $\varphi_{2k-1}(x)$ satisfying \eqref{1.200}, then $-\widetilde{a}_0$ is the coefficient on the $x$ term of
 $$
 g_{2k}(x)=1-x\varphi_{2k-1}(x)=\big(t^{2k}+t^{-2k}\big)^{-1}
 \cdot\mathcal{C}_{2k}\left(\frac{\kappa+1}{\kappa-1}-\frac{2x}{\alpha_1-\alpha_n}
\right),
$$
and $\widetilde{a}_0=-\frac{{\rm d}g_{2k}}{{\rm d}x} \big\vert_{x=0}$.  Combining \eqref{17.0003}, \eqref{eq1529}, \eqref{1932000}, Lemma \ref{lem19.59}, we complete the proof.
\end{proof}

We are ready to establish the following {\it nonasymptotic} bound on $\la_{opt}-\la_k$.
\begin{theorem}\label{8.42}
Suppose that $\|\bm x_{opt}\|=\|\bm x_k\|=\D$, $k=1,2,\ldots,k_{\max}$. Then
\begin{equation}\label{8.26}
0\leq\lambda_{opt}-\lambda_k\leq 2s_k\cdot\left(1+\frac{2k\sqrt{\kappa}}{\|A_{opt}\|}\left(\frac{\|\bm x_{opt}\|_{A_{opt}}}{\D}\right)^2\right)
\left(\frac{\sqrt{\kappa}-1}{\sqrt{\kappa}+1}\right)^{2k},
\end{equation}
where
\begin{equation}\label{eq5131}
s_k=\frac{\D^2}{\|\bm g\|^2\cdot\bm e_1^T(T_k+\la_{opt}I)^{-3}\bm e_1}.
\end{equation}
\end{theorem}
\begin{proof}

It follows from \eqref{9.44} and \eqref{19.56} that
\begin{align}
&\|\bm x_{k}\|^2-\|\widetilde{\bm x}_k\|^2=\|\bm g\|^2\cdot\bm e^T_1{\Big(}(T_{k}+\la_{k}I)^{-2} -(T_k+\la_{opt}I)^{-2}{\Big)}\bm e_1\nonumber\\
= &\|\bm g\|^2\!\cdot\!
\bm e^T_1{\Big(}(T_k\!+\!\la_kI)^{-1}\!-\!
(T_k\!+\!\la_{opt}I)^{-1}{\Big)}{\Big(}(T_k\!+\!\la_kI)^{-1}
\!+\!(T_k\!+\!\la_{opt}I)^{-1}{\Big)}\bm e_1\label{coro5.2}.
\end{align}
Notice that $\la_{k}\leq \la_{opt}$, so we have
\begin{eqnarray*}
\left\{
\begin{array}{ll}
(T_k+\la_kI)^{-1}-(T_k+\la_{opt}I)^{-1}=(\la_{opt}-\la_{k})(T_k+\la_kI)^{-1}(T_k+\la_{opt}I)^{-1},\\
\specialrule{0em}{0.8ex}{0.8ex}
(T_k+\la_kI)^{-1}(T_k+\la_{opt}I)^{-1}\big((T_k+\la_kI)^{-1}+(T_k+\la_{opt}I)^{-1}\big)
 \succeq 2(T_k+\la_{opt}I)^{-3}.
 \end{array}\right.
\end{eqnarray*}
Therefore,
\begin{align}
     &2~(\la_{opt}-\la_k)\cdot\bm e_1^T(T_k+\la_{opt}I)^{-3}\bm e_1\nonumber \\
\leq&\bm e^T_1{\Big(}(T_k+\la_kI)^{-1}-(T_k+\la_{opt}I)^{-1}{\Big)}{\Big(}(T_k+\la_kI)^{-1}
+(T_k+\la_{opt}I)^{-1}{\Big)}\bm e_1.\label{20.0000}
\end{align}
As $\|\bm x_{opt}\| =\|\bm x_{k}\|$, from \eqref{coro5.2}--\eqref{20.0000}, we obtain
\begin{align}\label{8.18}
0\!\overset{\eqref{eq9.42}}{\leq}\!\la_{opt}\!-\!\la_k\!\leq\!\frac{\|\bm x_{k}\|^2-\|\widetilde{\bm
x}_k\|^2}{2\|\bm g\|^2\cdot\bm e_1^T(T_k\!+\!\la_{opt}I)^{-3}\bm e_1}
=\frac{\|\bm x_{opt}\|^2-\|\widetilde{\bm
x}_k\|^2}{2\|\bm g\|^2\cdot\bm e_1^T(T_k\!+\!\la_{opt}I)^{-3}\bm e_1}.
\end{align}
A combination of Lemma \ref{22.22} yields \eqref{8.26}.
\end{proof}

Next, we focus on $s_k$ and prove that
\begin{equation}\label{eq513}
s_k=\frac{\D^2}{\|\bm g\|^2\cdot\bm e_1^T(T_k+\la_{opt}I)^{-3}\bm e_1}\rightarrow\frac{\D^2}{\bm g^TA^{-3}_{opt}\bm g}= s(\la_{opt})
\end{equation}
as $k$ increase,
moreover, we have that
$s_{k_{\max}}=s({\la}_{opt})$.
It is only necessary to consider the denominator $\|\bm g\|^2\cdot\bm e_1^T(T_k+\la_{opt}I)^{-3}\bm e_1$.
\begin{theorem}\label{thm20.08}
For $k=1,2,\ldots,k_{\max}$, we have
\begin{equation}\label{eq.2009}
0\!\leq\!\frac{\bm g^T\!A_{opt}^{-3}\bm g\!-\!\|\bm g\|^2\!\cdot\!
\bm e^T_1(T_{k}\!+\!\la_{opt}I)^{-3}\bm e_1}{\bm g^T\!A_{opt}^{-3}\bm g}\leq 4\left(1+\tau_{k}\right){\Big(}
\frac{\sqrt{{\kappa}}-1}{\sqrt{{\kappa}}+1}{\Big)}^{2k}
\end{equation}
and $\bm g^TA_{opt}^{-3}\bm g=\|\bm g\|^2\cdot\bm e^T_1(T_{k_{\max}}+\la_{opt}I)^{-3}\bm e_1$, where
\begin{align*}
\tau_{k}=\frac{\D^2(2k+1)\sqrt{\kappa}}{\bm g^T\!A_{opt}^{-3}\bm g\cdot\|A_{opt}\|}\left(1+\frac{2k\sqrt{\kappa}}{\|A_{opt}\|}\left(\frac{\|\bm x_{opt}\|_{A_{opt}}}{\D}\right)^2\right).
\end{align*}
Consequently, as $k$ increases, we have that
\begin{equation}
s_k\rightarrow s(\la_{opt})~~ and~~ s_{k_{\max}}=s(\la_{opt}).
\end{equation}
\end{theorem}
\begin{proof}
From $A_{opt}Q_{k_{\max}}\!=\!Q_{k_{\max}}(T_{k_{\max}}\!+\!\la_{opt}I)$, we have $A_{opt}^{-3}Q_{k_{\max}}\!=\!Q_{k_{\max}}(T_{k_{\max}}\!+\!\la_{opt}I)^{-3}$. As $Q_{k_{\max}}\bm e_1=\frac{\bm g}{\|\bm g\|}$, we have from \cite[Theorem 4.1]{GG.} that
\begin{align}\label{eq2027}
\bm g^T\!A_{opt}^{-3}\bm g\!-\!\|\bm g\|^2\!\cdot\!
\bm e^T_1T_{k,opt}^{-3}\bm e_1\!=\!\|\bm g\|^2\cdot\left(\bm e_1^T(T_{k_{\max}}\!+\!\la_{opt}I)^{-3}\bm e_1\!-\!\bm
e_1^TT_{k,opt}^{-3}\bm e_1\right)\!\geq \!0,
\end{align}
and $\bm g^TA_{opt}^{-3}\bm g=\|\bm g\|^2\cdot\bm e^T_1(T_{k_{\max}}+\la_{opt}I)^{-3}\bm e_1$, i.e., $s(\la_{k_{\max}})=s(\la_{opt})$,
where $T_{k,opt}=T_k+\la_{opt}I$.

Suppose that $\varphi_{2k}(x)\in \mathcal{P}_{2k}$ satisfies \eqref{1.200}. Let
$$
q_{2k-2}(x)=\frac{\varphi_{2k}(x)-\widetilde{a}_1x-\widetilde{a}_0}{x^2}\in\mathcal{P}_{2k-2},~~{\rm with}~~x\in[\alpha_n+\la_{opt},\alpha_1+\la_{opt}],
$$
where $\widetilde{a}_0$ and $\widetilde{a}_1$ are the constant term and the coefficient on the $x$ term of $ \varphi_{2k}(x)$, respectively.
By Lemma \ref{lem14.43},
\begin{small}
\begin{align}
    \frac{\bm g^TA_{opt}^{-3}\bm g}{\|\bm g\|^2}-\bm e^T_1T_{k,opt}^{-3}\bm e_1=\frac{\bm g^T\big[A_{opt}^{-3}-q_{2k-2}(A_{opt})\big]\bm g}{\|\bm g\|^2}-\bm e^T_1\big[T_{k,opt}^{-3}-q_{2k-2}(T_{k,opt})\big]\bm e_1.\label{eq15.22}
\end{align}
\end{small}

Now we consider $\bm g^T\big[A_{opt}^{-3}-q_{2k-2}(A_{opt})\big]\bm g$. From the definition of $q_{2k-2}(x)$, we obtain
\begin{align}
    &\bm g^T\big[A_{opt}^{-3}-q_{2k-2}(A_{opt})\big]\bm g
=\bm g^TA^{-\frac{3}{2}}_{opt}\big[I-A_{opt}^3q_{2k-2}(A_{opt})
\big]A^{-\frac{3}{2}}_{opt}\bm g\nonumber\\
=&\bm g^TA^{-\frac{3}{2}}_{opt}
\big[I-A_{opt}^3q_{2k-2}(A_{opt})-\widetilde{a}_1A^2_{opt}
-\widetilde{a}_0A_{opt}\big]A^{-\frac{3}{2}}_{opt}\bm g+\widetilde{a}_1\bm g^TA^{-1}_{opt}\bm g+\widetilde{a}_0\bm g^TA^{-2}_{opt}\bm g\nonumber\\
=&\bm g^TA^{-\frac{3}{2}}_{opt}\big[I-A_{opt}\varphi_{2k}(A_{opt})\big]
A^{-\frac{3}{2}}_{opt}\bm g-\widetilde{a}_1\bm g^T\bm x_{opt}+\widetilde{a}_0\|\bm x_{opt}\|^2.
\end{align}
Similarly, let ${\widetilde{\bm h}}_k$ be defined in \eqref{19.56}, we arrive at
\begin{align}
&\bm e_1^T\!\big[T_{k,opt}^{-3}-q_{2k-2}(T_{k,opt})\big]\bm e_1\nonumber\\
=&\bm e_1^TT_{k,opt}^{-\frac{3}{2}}\big[I\!-\!
T_{k,opt}\varphi_{2k}(T_{k,opt})\big]
T_{k,opt}^{-\frac{3}{2}}\bm e_1-\widetilde{a}_1\frac{\bm e_1^T\widetilde{\bm h}_{k}}{\|\bm g\|}\!+\!\widetilde{a}_0\frac{\| \widetilde{\bm h}_{k}\|^2}{\|\bm g\|^2}.\label{eq1604}
\end{align}
A combination of \eqref{eq15.22}--\eqref{eq1604} yields
\begin{small}
\begin{align}
    &\bm g^TA_{opt}^{-3}\bm g-\|\bm g\|^2\cdot\bm e^T_1(T_{k}+\la_{opt}I)^{-3}\bm e_1\nonumber\\
  =&\bm g^TA^{-\frac{3}{2}}_{opt}\big[I\!-\!A_{opt}\varphi_{2k}(A_{opt})\big]
  A^{-\frac{3}{2}}_{opt}\bm g\!-\!\|\bm g\|^2\bm e_1^TT_{k,opt}^{-\frac{3}{2}}\big[I\!-\!T_{k,opt}\varphi_{2k}(T_{k,opt})\big]
T_{k,opt}^{-\frac{3}{2}}\bm e_1\nonumber\\
   &~~+\widetilde{a}_1(\|\bm g\|\cdot\bm e_1^T\widetilde{\bm h}_{k}-\bm g^T\bm x_{opt})+\widetilde{a}_0(\|\bm x_{opt}\|^2-\|{\widetilde{\bm h}}_k\|^2)\nonumber\\
   \leq&\bm g^TA^{-3}_{opt}\bm g\cdot \|I-A_{opt}\varphi_{2k}(A_{opt})\|+\bm e_1^T T_{k,opt}^{-{3}}\bm e_1\cdot\|\bm g\|^2\|I-T_{k,opt}\varphi_{2k}(T_{k,opt})\|
\nonumber\\
   &~~+\widetilde{a}_1(\|\bm g\|\cdot\bm e_1^T\widetilde{\bm h}_{k}-\bm g^T\bm x_{opt})+\widetilde{a}_0(\|\bm x_{opt}\|^2-\|{\widetilde{\bm h}}_k\|^2)\nonumber\\
   \overset{\eqref{eq2027}}{\leq}&\bm g^TA^{-3}_{opt}\bm g\cdot\big( \|I-A_{opt}\varphi_{2k}(A_{opt})\|+\|I-T_{k,opt}\varphi_{2k}(
   T_{k,opt})\|\big)
\nonumber\\
   &~~+\widetilde{a}_1(\|\bm g\|\cdot\bm e_1^T\widetilde{\bm h}_{k}-\bm g^T\bm x_{opt})+\widetilde{a}_0(\|\bm x_{opt}\|^2-\|{\widetilde{\bm h}}_k\|^2)\nonumber\\
  \overset{\rm \eqref{eq15,35}}{\leq}&4\bm g^TA^{-3}_{opt}\bm g\!\cdot\!{\Big(}\frac{\sqrt{\kappa}-1}
   {\sqrt{\kappa}+1}{\Big)}^{2k}\!+\!\widetilde{a}_1(\|\bm g\|\!\cdot\!\bm e_1^T\widetilde{\bm h}_{k}\!-\!\bm g^T\!\bm x_{opt})\!+\!\widetilde{a}_0(\|\bm x_{opt}\|^2\!-\!\|{\widetilde{\bm x}}_k\|^2),\label{eq15.29}
\end{align}
\end{small}
where we used
$\|{\widetilde{\bm h}}_k\|\!=\!\|{\widetilde{\bm x}}_k\|$; see \eqref{19.56}.

In terms of the definition of $\varphi_{2k}(x)$, we see that $-\widetilde{a}_0$ and $-\widetilde{a}_1$ are the coefficients on the terms $x$ and $x^2$ of
$$
\widetilde{g}_{2k+1}(x)=1-x\varphi_{2k}(x)=
\big(t^{2k+1}+t^{-(2k+1)}\big)^{-1}\cdot
\mathcal{C}_{2k+1}\left(\frac{\kappa+1}{\kappa-1}-\frac{2x}{\alpha_1-\alpha_n} \right)
,
$$
respectively. Thus,
$$
\widetilde{a}_0=-\frac{{\rm d}\widetilde{g}_{2k+1}}{{\rm d}x} \bigg\vert_{x=0}~~{\rm and}~~\widetilde{a}_1=-\frac{1}{2}~\frac{{\rm d}^2\widetilde{g}_{2k+1}}{{\rm d}x^2} \bigg\vert_{x=0}.$$
Moreover, it follows from Lemma \ref{lem19.59} and \eqref{1932001} that
\begin{equation}\label{eq532}
0\!\leq\!\widetilde{a}_0\!\leq\!\frac{(2k\!+\!1)\sqrt{\kappa}}{\|A_{opt}\|}
~~{\rm and}~~\widetilde{a}_1\big(\|\bm g\|\cdot\bm e_1^T\widetilde{\bm h}_{k}-\bm g^T\bm x_{opt}\big)\!\leq\!0.
\end{equation}
A combination of \eqref{eq15.29}, \eqref{eq532}, and Lemma \ref{22.22} yields \eqref{eq.2009}.
\end{proof}
\begin{remark}\label{remk16.59}
In \cite[eq. (4.9)]{5}, the condition number of $\la_{opt}$ is defined as
$$
cond(\la_{opt})=\frac{1}{2}\vert\bm y_1^T\bm y_2\vert^{-1}.
$$
As $\bm x_{opt}=-sign(\bm g^T\bm y_2)\frac{\D\cdot\bm y_1}{\|\bm y_1\|}$ \cite[eq. (21)]{3} and $M\bm y=\la_{opt}\bm y$, we have $\bm y_1=A_{opt}\bm y_2$ and
\begin{align}
{s}(\la_{opt})\!=\!\frac{\D^2}{\bm g^T\!A^{-3}_{opt}\bm g}=\frac{\D^2}{\bm x_{opt}^TA^{-1}_{opt}\bm x_{opt}}
\!=\!\frac{\|\bm y_1\|^2}{\bm y_1^T\!A_{opt}^{-1}\bm y_1}
\!=\!\frac{\|\bm y_1\|^2}{\vert\bm y_1^T\!\bm y_2\vert}
\!=\!2\|\bm y_1\|^2 cond(\la_{opt}).\label{eq2009}
\end{align}
Thus, $s(\la_{opt})$ is closely related to the condition number $cond(\la_{opt})$.
However, we point out that $s(\la_{opt})$ can be much smaller than $cond(\la_{opt})$ in practice. This can happen, say, when
$\|\bm y_1\|\ll 1$ in the ``nearly hard cases'' \cite{3}; see the numerical results in Table \ref{tab:16.39}.
That is, $s(\la_{opt})$ can be much different from the condition number $cond(\la_{opt})$, and it is more appropriate for depicting the convergence of $\lambda_k$.
\end{remark}


\begin{remark}\label{14.55}
{\it In summary, Theorem \ref{8.42} and Theorem \ref{thm20.08} show the importance of $s(\lambda_{opt})$ on the convergence of $\lambda_k$.
More precisely, $\lambda_k$ may converge slow if $s(\lambda_{opt})$ is very large; see Figure \ref{fig19.1201}.
Compared with the result given in Theorem \ref{14.144}, our
 bound \eqref{8.26} is {\it non-asymptotic}, and there is no need to assume that  $k$ is sufficiently large.}
\end{remark}

\section{On the convergence of Krylov subspace method for the cubic regularization problem}
In this section, we consider the following cubic regularization problem, which is a regularization variant of the TRS \eqref{1} \cite{T.B,0,110,G.C.P,C.G.M,G.S.,J&Z,F.L,Y.N}:
\begin{equation}\label{2}
\min_{\bm x\in \mathbb{R}^n} f_{\sigma}(\bm x)~~{\rm with}~~f_{\sigma}(\bm x)= \frac{1}{2}\bm x^TA\bm x+\bm x^T\bm g+\frac{\sigma}{3}\|\bm x\|^3,
\end{equation}
where $\sigma>0$.
Denote by $\mu_{opt}=\sg\|\bm x^{\sg}_{opt}\|$, and by
\begin{equation}\label{eq.1458}
\e_{\sg,k}=\|(I-Q_kQ_k^T)\bm x_{opt}^{\sg}\|,~A_{\sg,opt}=A+\mu_{opt}I~~{\rm and}~~\kappa_{\sg}=\|A_{\sg,opt}\|\|A^{-1}_{\sg,opt}\|.
\end{equation}
We have
\begin{equation}\label{eqsg}
\e_{\sg,k}=\|(I-Q_kQ^T_k)\bm x^{\sg}_{opt}\|\leq 2\|\bm x^{\sg}_{opt}\|
\cdot {\bigg(}\frac{\sqrt{\kappa_{\sg}}-1}{\sqrt{\kappa_{\sg}}+1}{\bigg)}^{k},
\end{equation}
whose proof is similar to that of \eqref{11.09}.
We omit the details and refer to \cite{5,4}.

The following theorem provides a necessary and sufficient condition for determining a global
optimal solution of the subproblem \eqref{2}.
\begin{theorem}\cite[Theorem 3.1]{G.C.P}\label{Thm61}
Any $\bm x_{opt}^{\sg}$ is a global minimizer of $f_{\sg}(\bm x)$ over $\mathbb{R}^n$ if and only if it satisfies
the system of equation
\begin{align}
(A+\mu_{opt}I){\bm x^{\sg}_{opt}}=-\bm g~~{ and}~~A+\mu_{opt}I\succcurlyeq \bm O\nonumber.
\end{align}
If $A+\mu_{opt}I$ is positive definite, $\bm x^{\sg}_{opt}$ is unique.
\end{theorem}

Suppose that $A+\mu_{opt}I\succ\bm O$, then according to Theorem \ref{Thm61}, $\bm x^{\sg}_{opt}$ is unique and the cubic regularization problem \eqref{2} is in \emph{easy case} \cite[Definition 2.3]{J&Z}. Otherwise, one fails to get a good approximation to $\bm x_{opt}^{\sg}$ from the Krylov subspace $\mathcal{K}_k(A,\bm g)$ \cite{G.C.P,G.S.}. In the Krylov subspace method, we solve the following problem \cite{J&Z}
 \begin{equation}\label{9.0000}
 \min_{\bm x\in \mathcal{K}_k(A,\bm g)}\left\{f_{\sg}(\bm x)=\frac{1}{2}\bm x^T A\bm x+\bm x^T\bm g+\frac{\sg}{3}\|\bm x\|^3\right\}
 \end{equation}
for the cubic regularized quadratic subproblem.

Let $Q_k$ be an orthonormal
basis for the Krylov subspace $\mathcal{K}_k(A,\bm g)$, and let $T_k$ be the tridiagonal matrix obtained from the $k$-step Lanczos process; refer to \eqref{21.18}. Let
 \begin{equation}\label{15.18}
 \bm h_k^{\sg}=\arg\min_{\bm h\in\mathbb{R}^k}\left\{\frac{1}{2}\bm h^TT_k\bm h+\|\bm g\|\bm h^T\bm e_1+\frac{\sg}{3}\|\bm h\|^3\right\},
 \end{equation}
which satisfies
 \begin{equation}\label{eq1640}
(T_k+\sg \|\bm h_{k}^{\sg}\|I){\bm h^{\sg}_{k}}=-\|\bm g\|\bm e_1~~{\rm and}~~
T_k+\sg \|\bm h_k^{\sg}\|I \succcurlyeq \bm O.
\end{equation}
As $T_k$ is a irreducible tridiagonal matrix, by \cite[Theorem 5.3]{10} and \cite[Definition 2.3]{J&Z}, the cubic regularized quadratic subproblem \eqref{15.18} is also in $\emph{easy case}$, where $T_k+\sg \|\bm h_k^{\sg}\|I \succ \bm O$ and $\bm h^\sg_{k}=-\|\bm g\|\cdot(T_k+\sg \|\bm h_k^{\sg}\|I)^{-1}\bm e_1$ is unique.
Denote by $\bm x^{\sg}_{k}=Q_k\bm h_k^{\sg}$, then it can be utilized as an approximation to $\bm x^{\sg}_{opt}$.
By \cite[Theorem 1]{C.G.M}, we see that
\begin{equation}\label{eq.19.45}
\|\bm x^\sg_1\|\leq\|\bm x^\sg_2\|\leq\cdots\leq\|\bm x^\sg_{k_{\max}}\|=\|\bm x^\sg_{opt}\|.
\end{equation}
Similar to \eqref{11.10},
we have
 \begin{equation}\label{11.100}
\left\{
\begin{array}{ll}
\bm x^{\sg}_k \in  \mathcal{K}_k(A,\bm g)                  \\
\bm r^{\sg}_k=(A+\sg \|\bm x^{\sg}_k\|I)\bm x^{\sg}_k+\bm g~ \bot ~\mathcal{K}_k(A,\bm g).
\end{array}\right.
\end{equation}

In this section, we are interested in the convergence of the Krylov subspace method for the subproblem \eqref{2}. Some upper bounds on $\|(A+\mu_kI)\bm x^{\sg}_{k}+\bm g\|$ and $\vert\mu_k-\mu_{opt}\vert$ were established in \cite{G.S.,J&Z}, where $\mu_k=\sg \|\bm x^{\sg}_k\|$. We try to establish some sharper error bounds in this section.
\subsection{Some existing results}
The following result is due to Gould and Simoncini \cite{G.S.}:
\begin{theorem}\cite[Theorem 3.4]{G.S.}
The residual $(A+\mu_{k}I)\bm x^\sg_k+\bm g$ for the k-th iterate, $\bm x^\sg_k$, generated by the regularization subproblem \eqref{9.0000} satisfies the
bound
\begin{equation}\label{20.32}
\|(A+\mu_{k}I)\bm x^{\sg}_k+\bm g\| \leq\|\bm g\|{\Big(}\frac{2\beta_k\kappa_{\sg,k}}{\|T_k+\mu_kI\|}{\Big)} {\Big(}\frac{\sqrt{\kappa_{\sg,k}}-1}{\sqrt{\kappa_{\sg,k}}+1}{\Big)}^{k-1},
\end{equation}
where $\kappa_{\sg,k}$ is the 2-condition number of $T_k+\mu_kI$ and $\beta_k$ is the $(k,k+1)$-st entry of $T_{k+1}$.
\end{theorem}

 Recently, X. Jia {\it et al.} \cite{J&Z} established the following bound on $\vert\mu_k^3-\mu_{opt}^3\vert$:
 \begin{theorem}\label{thm15.12}\cite[Theorem 4.6]{J&Z}
Suppose that the cubic regularization problem \eqref{2} is in \emph{easy case}. Then
 \begin{equation}\label{eq2054}
 |\mu_k^3-\mu_{opt}^3|\leq6\sg^2\|\bm g\|\sqrt{\kappa}\cdot\Pi_k+12\sg^2\|A_{\sg,opt}\|\cdot\Pi_k^2,
\end{equation}
where
\begin{equation*}
\Pi_k=\frac{2\|\bm g\|(\varsigma_\sg+\sqrt{\varsigma_\sg^2-1})^{-k}}{(\alpha_1-\alpha_n)(\varsigma_\sg^2-1)}
~~{ with}~~\varsigma_\sg=\frac{\kappa_{\sg}+1}{\kappa_{\sg}-1}.
\end{equation*}
 \end{theorem}
\subsection{Improved upper bounds on $\bm{\|(A+ \mu_{k} I) x_{\sg,k}+g\|}$ and $\bm{\mu_{opt}-\mu_k}$}
We are in a position to establish a new upper bound on $\|(A+\mu_{k}I)\bm x^{\sg}_k+\bm g\|$.
\begin{theorem}
Suppose that the cubic regularization problem \eqref{2} is in \emph{easy case}. Then
\begin{equation}\label{20.41}
\|(A+\mu_kI)\bm x^{\sg}_{k}+\bm g\|\leq \min\{\xi_{\sg 1},\xi_{\sg 2}\},
\end{equation}
where
$$
\xi_{\sg 1}\!=\!2\|A_{\sg,opt}\|\sqrt{\|\bm x_{opt}^{\sg}\|^2\!+\!\e_{\sg,k}^2}\cdot{\Big(}\frac{\sqrt{\kappa_{\sg}}\!
-\!1}{\sqrt{\kappa_{\sg}}\!+\!1}{\Big)}^{k}
~~ and~~
\xi_{\sg 2}\!=\!2\beta_k\|\bm x_{k}^{\sg}\|\cdot{\Big(}\frac{\sqrt{\kappa_{\sg,k}}\!-\!1
}{\sqrt{\kappa_{\sg,k}}\!+\!1}{\Big)}^{k-1},
$$
with $\beta_k$ being the $(k,k+1)$-st entry of $T_{k+1}$ and $\kappa_{\sg,k}$ being the 2-condition number of $T_k+\mu_kI$.
\end{theorem}
\begin{proof}
Without loss of generality, we suppose that $\bm x_{opt}^{\sg}\neq\bm 0$. Otherwise, we have from \eqref{eq.19.45} that $\bm x_{k}^{\sg}=\bm 0$.
Denote by
$$
\rho(\bm x)=\frac{\mu_{opt}}{2}\cdot(\|\bm x^{\sg}_{opt}\|^2-\|\bm x\|^2)+\frac{\sg}{3}\cdot(\|\bm x\|^3-\|\bm x^{\sg}_{opt}\|^3).
$$
Recall that $\mu_{opt}=\sg\|\bm x_{opt}^\sg\|$. Then
\begin{align}
\rho(\bm x_k^\sg)=&\frac{\sg\|\bm x^{\sg}_{opt}\|}{2}(\|\bm x^{\sg}_{opt}\|^2-\|\bm x_k^{\sg}\|^2)+\frac{\sg}{3}(\|\bm x_k^{\sg}\|^3-\|\bm x^{\sg}_{opt}\|^3)\nonumber\\
=&\frac{\sg}{6}{\big (}\|\bm x^{\sigma}_{opt}\|^2+\|\bm x^{\sigma}_{opt}\|\|\bm x^{\sigma}_{k}\|-2\|\bm x^{\sigma}_{k}\|^2{\big )}(\|\bm x^{\sigma}_{opt}\|
-\|\bm x^{\sigma}_{k}\|)\overset{\eqref{eq.19.45}}{\geq} 0.\label{eq1548}
\end{align}
 For any $\bm x$, we have that
\begin{small}
\begin{align}
0&\leq f_{\sg}(\bm x)- f_{\sigma}(\bm x^{\sigma}_{opt})=
\frac{1}{2}\left(\bm x^TA\bm x-(\bm x^{\sigma}_{opt})^TA\bm x^{\sigma}_{opt}\right)
+\bm g^T(\bm x-\bm x^{\sg}_{opt})+\frac{\sg}{3}(\|\bm x\|^3-\|\bm x^{\sg}_{opt}\|^3)\nonumber\\
         &=\frac{1}{2}\left(\bm x^TA_{\sg,opt}\bm x-(\bm x^\sg_{opt})^TA_{\sg,opt}\bm x^\sg_{opt}\right)
         -(\bm x^\sg_{opt})^TA_{\sg,opt}(\bm x-\bm x^\sg_{opt})+\rho(\bm x)\nonumber\\
&=\frac{1}{2}(\bm x-\bm x^{\sg}_{opt})^TA_{\sg, opt}(\bm x-\bm x^{\sg}_{opt})+\rho(\bm x).\label{eq1056}
\end{align}
\end{small}

Let $\bm y^{\sg}_{k}\!=\!\frac{Q_kQ_k^T\!\bm x^{\sg}_{opt}}{\|Q_kQ_k^T\!\bm x^{\sg}_{opt}\|}\|\bm x^{\sg}_{opt}\|$. As $\vert\bm g^T\!\bm x^{\sg}_{opt}\vert=(\bm
x^{\sg}_{opt})^TA_{\sg,opt}\bm x^{\sg}_{opt}\!>\!0$, we have $Q_kQ^T_k\bm x^{\sg}_{opt}\neq\bm 0$, and $\bm y^{\sg}_k$ is well-defined. Notice that $\rho(\bm y^{\sg}_{k})=0$, so we have
\begin{align}
0\!&\leq\!f_{\sg}(\bm x^{\sg}_{k})-f_{\sigma}(\bm x^{\sigma}_{opt})\leq
f_{\sg}(\bm y^{\sg}_{k})\!-\!f_{\sigma}(\bm x^{\sigma}_{opt})
\!\leq\!\frac{1}{2}\|A_{\sg, opt}\|\| \bm
y^{\sg}_{k}\!-\!\bm x^{\sg}_{opt}\|^2\nonumber\\
&\overset{\eqref{10.24}-\eqref{eq10.29}}{\leq} \frac{1}{2}\|A_{\sg,opt}\|{\Big (}1\!+\!\frac{\epsilon_{\sg,k}^2}{\|\bm x^{\sigma}_{opt}\|^2} {\Big )}
\epsilon_{\sg,k}^2\nonumber\\
&\overset{\eqref{eqsg}}{\leq} 2\|A_{\sg,opt}\|{\Big(}\|\bm x^{\sigma}_{opt}\|^2\!+\!{\epsilon_{\sg,k}^2}{\Big)}
{\Big(}\frac{\sqrt{\kappa_\sg}\!-\!1}{\sqrt{\kappa_\sg}\!+\!
1}{\Big)}^{2k}.\label{eq1054}
\end{align}

From \eqref{11.59} and \eqref{11.100}, we obtain
\begin{align}
\|(A\!+\!\mu_{k}I)\bm x^{\sg}_k\!+\!\bm g\|^2\!=&\|A_{\sg,opt}(\bm x^{\sg}_k-\bm x^{\sg}_{opt})\|^2-(\mu_{opt}-\mu_{k})^2\|\bm x^{\sg}_{opt}\|^2\nonumber\\
\leq&\|A_{\sg,opt}(\bm x^{\sg}_k\!-\!\bm x^{\sg}_{opt})\|^2\!\leq\!
\|A_{\sg,opt}
\|\!\cdot\!(\bm x_k^{\sg}\!-\!\bm x^{\sg}_{opt})^T\!A_{\sg, opt}(\bm x_k^{\sg}\!-\!\bm x^{\sg}_{opt})\nonumber\\
\overset{\eqref{eq1548}}{\leq} &2\|A_{\sg,opt}\|\cdot\Big(\frac{1}{2}(\bm x_k^{\sg}-\bm x^{\sg}_{opt})^TA_{\sg, opt}(\bm x_k^{\sg}-\bm
x^{\sg}_{opt})+\rho(\bm x_k^\sg)\Big)\nonumber\\
\overset{\eqref{eq1056}}{=}&2\|A_{\sg,opt}\|\cdot\left( f_{\sg}(\bm x^{\sg}_{k})\!-\!f_{\sigma}(\bm x^{\sigma}_{opt})\right)\overset{\eqref{eq1054}}{\leq} \xi_{\sg 1}^2.
\end{align}

  Recall that $T_k\!+\!\mu_k I\!\succ\!\bm O$ is a symmetric and tridiagonal matrix, and $\bm h_k^\sg\!=\!-\|\bm g\|(T_k\!+\!\mu_k I)^{-1}\bm e_1$. Similar to the proof of Theorem \ref{9.051}, we have $\|(A\!+\!\mu_{k}I)\bm x^{\sg}_k\!+\!\bm g\|\!\leq \!\xi_{\sg 2}$.
This completes the proof.
\end{proof}
\begin{remark}
{\it It follows from \eqref{eq1640} that
$$
\|\bm x^{\sg}_k\|=\|\bm h^{\sg}_{k}\|=\|\bm g\|\|(T_k+\mu_kI)^{-1}\bm e_1\|\leq
{\|(T_k+\mu_kI)^{-1}\|\|\bm g\|}=\frac{\|\bm g\|\kappa_{\sg,k}}{\|T_k+\mu_kI\|},
$$
so our bound \eqref{20.41} is no worse than \eqref{20.32}. Moreover, in  \eqref{20.41}, there is no $\kappa_{\sg,k}$ in the coefficient before ${\Big(}\frac{\sqrt{\kappa_{\sg}}\!
-\!1}{\sqrt{\kappa_{\sg}}\!+\!1}{\Big)}^k$  and ${\Big(}\frac{\sqrt{\kappa_{\sg,k}}\!-\!1
}{\sqrt{\kappa_{\sg,k}}\!+\!1}{\Big)}^k$.
}
\end{remark}

Next, we consider the convergence of $\mu_k$.
\begin{theorem}\label{thm20.06}
Suppose that the cubic regularization problem \eqref{2} is in \emph{easy case}. Then
\begin{equation*}
0\leq\mu_{opt}-\mu_k\leq 2 s_{\sg,k}\cdot\left(1+\frac{2k\sqrt{\kappa_{\sg}}}{\|A_{\sg,opt}\|}\left(\frac{\|\bm x^\sg_{opt}\|_{A_{\sg,opt}}}{\|\bm x_{opt}^\sg\|}\right)^2\right)
\left(\frac{\sqrt{\kappa_{\sg}}-1}{\sqrt{\kappa_{\sg}}+1}\right)^{2k},
\end{equation*}
where
$
s_{\sg,k}=\frac{\|\bm x_{opt}^\sg\|^2}{\|\bm g\|^2\cdot\bm e_1^T(T_k+\mu_{opt}I)^{-3}\bm e_1}.
$
\end{theorem}
\begin{proof}
Recall that $\mu_{opt}=\sg \|\bm x_{opt}^{\sg}\|$, $\mu_k=\sg \|\bm x_{k}^{\sg}\|$ and $\|\bm x_{opt}^{\sg}\|\geq\|\bm x_{k}^{\sg}\|$. Hence,
$\mu_{opt}\geq \mu_k$ and $T_k+\mu_{opt}I\succcurlyeq T_k+\mu_{k}I\succ\bm O$.
Let $\widetilde{\bm x}^{\sg}_{k}=-\|\bm g\|\cdot Q_k(T_k+\mu_{opt}I)^{-1}\bm e_1$.
From \eqref{coro5.2}--\eqref{8.18}, we have
$$
0\leq\mu_{opt}-\mu_k\leq\frac{\|\bm x^{\sg}_{k}\|^2-\|\widetilde{\bm x}^{\sg}_k\|^2}{2\|\bm g\|^2\cdot\bm e_1^T(T_k+\mu_{opt}I)^{-3}\bm e_1}\overset{\eqref{eq.19.45}}{\leq}\frac{\|\bm
x^{\sg}_{opt}\|^2-\|\widetilde{\bm x}^{\sg}_k\|^2}{2\|\bm g\|^2\cdot\bm e_1^T(T_k+\mu_{opt}I)^{-3}\bm e_1}.
$$
The remaining proof is similar to that of Theorem \ref{8.42}.
\end{proof}
\begin{corollary}\label{Cor69}
Suppose that the cubic regularization problem \eqref{2} is in \emph{easy case}. Then
\begin{equation}\label{eq.8.26}
0\leq\mu^3_{opt}-\mu^3_k\leq 6s_{\sg,k}\cdot {\Bigg(}\sg^2\|\bm x_{opt}^\sg\|^2+\frac{2\sg^2k\sqrt{\kappa_{\sg}}\|\bm x_{opt}^\sg\|^2_{A_{\sg,opt}}}{\|A_{\sg,opt}\|}{\Bigg)}
\left(\frac{\sqrt{\kappa_{\sg}}-1}{\sqrt{\kappa_{\sg}}+1}\right)^{2k},
\end{equation}
where
$
s_{\sg,k}=\frac{\|\bm x_{opt}^\sg\|^2}{\|\bm g\|^2\cdot\bm e_1^T(T_k+\mu_{opt}I)^{-3}\bm e_1}.
$
\end{corollary}
\begin{proof}
{As $\mu_{opt}\geq\mu_k\geq 0$ and $\mu_{opt}=\sg\|\bm x^\sg_{opt}\|$}, we have that
$$
\mu_{opt}^3\!-\!\mu_{k}^3\!=\!(\mu_{opt}\!-\!\mu_{k})(\mu_{opt}^2\!+\!\mu_{k}^2
\!+\!\mu_{opt}\mu_{k})\!\leq\!3\mu^2_{opt}(\mu_{opt}\!-\!\mu_{k})\!=\!3\sg^2\|\bm x^\sg_{opt}\|^2\!\cdot\!(\mu_{opt}\!-\!\mu_{k}),
$$
and a combination of Theorem \ref{thm20.06} yields the result.
\end{proof}
\begin{remark}
{\it Theorem \ref{thm15.12} shows that the upper bound of $\vert\mu_k^3-\mu_{opt}^3\vert$ is in the order of
$\Pi_k=\mathcal{O}\Big({\big(}\frac{\sqrt{\kappa_{\sg}}-1}
{\sqrt{\kappa_{\sg}}+1}{\big)}^{k}\Big)$ \cite{J&Z}.
As a comparison, Corollary \ref{Cor69} indicates that it is
in the order of $\mathcal{O}\Big({\big(}\frac{\sqrt{\kappa_{\sg}}-1}{\sqrt{\kappa_{\sg}}+1}{\big)}^{2k}\Big)$. Therefore, our new bound is much better than that of X. Jia et al.; see Figure \ref{fig:2124}.}
\end{remark}

\section{Numerical experiments}
In this section, we perform some numerical experiments to illustrate the effectiveness of our theoretical results.
 All the numerical experiments were run on a Inter(R) Core i5 with CPU 2.60Hz and RAM 4GB under Windows 7 operation system.
 The experimental results are obtained from using MATLAB R2014b implementation with machine precision ${\bf u}\approx 2.22\times 10^{-16}$.

{\bf Example 1.}
In this example, we compare our new bounds on $\|(A+\la_k I)\bm x_k+\bm g\|$, $\sin\angle(\bm x_{opt}, \bm x_{k})$, $\la_{opt}-\la_k$
with those of Z. Jia {\it et al}. \cite{5} and Gould {\it et al.} \cite{G.S.}.
We consider the diagonal matrix $A=diag(t^{[a,b]}_{jn})$ \cite{5,4}, with
\begin{equation}\label{eqn61}
t_{jn}^{[a,b]}={\Big(} \frac{b-a}{2}  {\Big )} {\Big(} t_{jn}+{\Big(} \frac{a+b}{b-a} {\Big )} {\Big )}
~~
{\rm with}
~~
t_{jn} = \cos \frac{(2j-1)\pi}{2n},~j=1,2,\ldots,n,
\end{equation}
where $n=10000$ and $t_{jn}^{[a,b]}$ is the $n$-th translated Chebyshev zero nodes on $[a,b]$. The vector $\bm g$ is generated by the MATLAB build-in function
${\tt randn.m}$, and is normalized with its Euclidean norm.

\begin{figure}[ht]\caption{\footnotesize Example 1: A comparison of the upper bounds on $\|(A+\la_k I)\bm x_k+\bm g\|$.}\label{fig:1449}
\centering
\subfigure[\footnotesize  $A=diag(t^{[-5,5]}_{jn})$, $\D
=1$.]{\includegraphics[width=6.3cm,height=3.8cm]{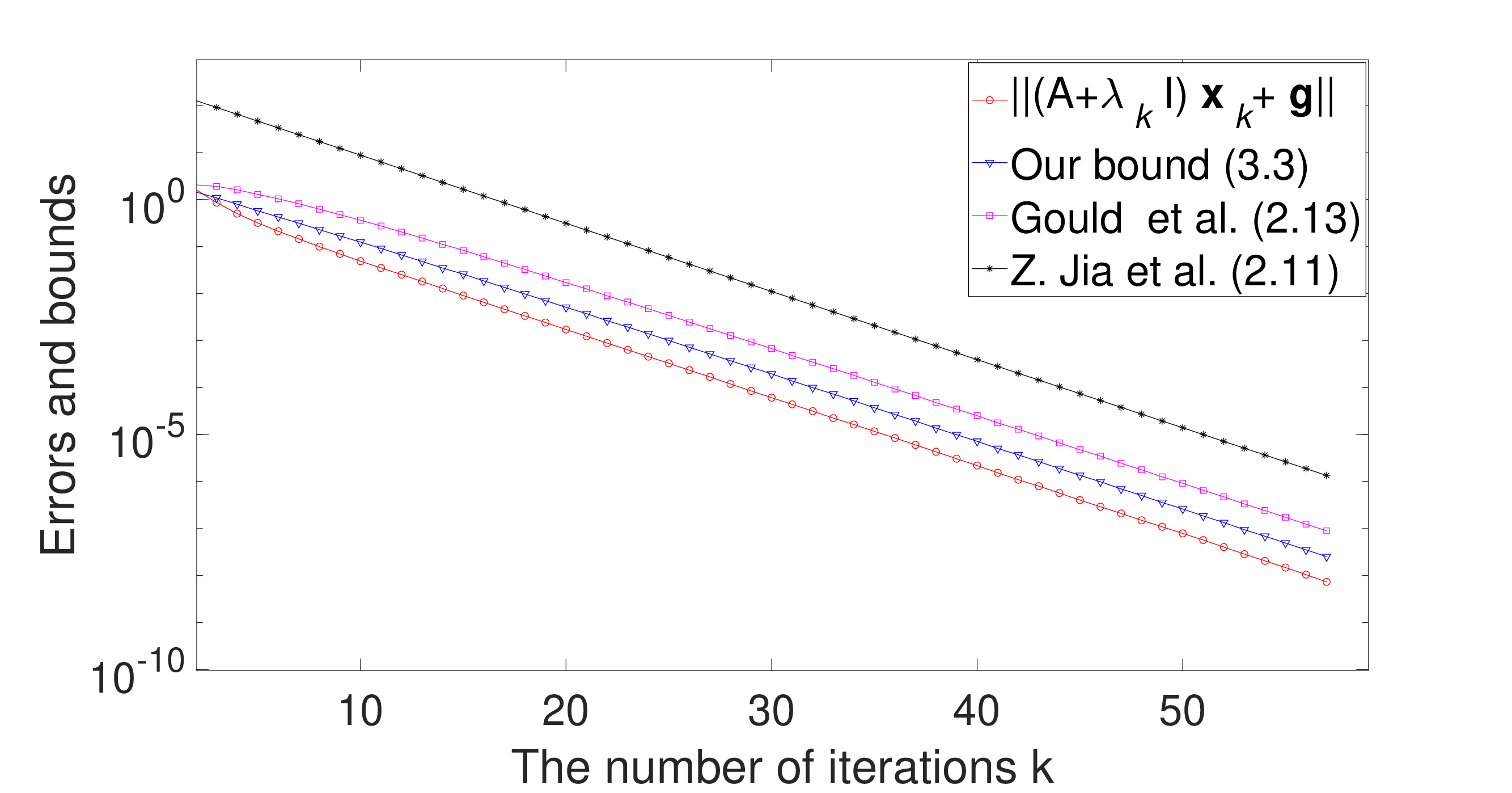}}~~~~
\subfigure[\footnotesize  $A=diag(t^{[-10,10]}_{jn})$, $\D=10$.]{\includegraphics[width=6.3cm,height=3.8cm]{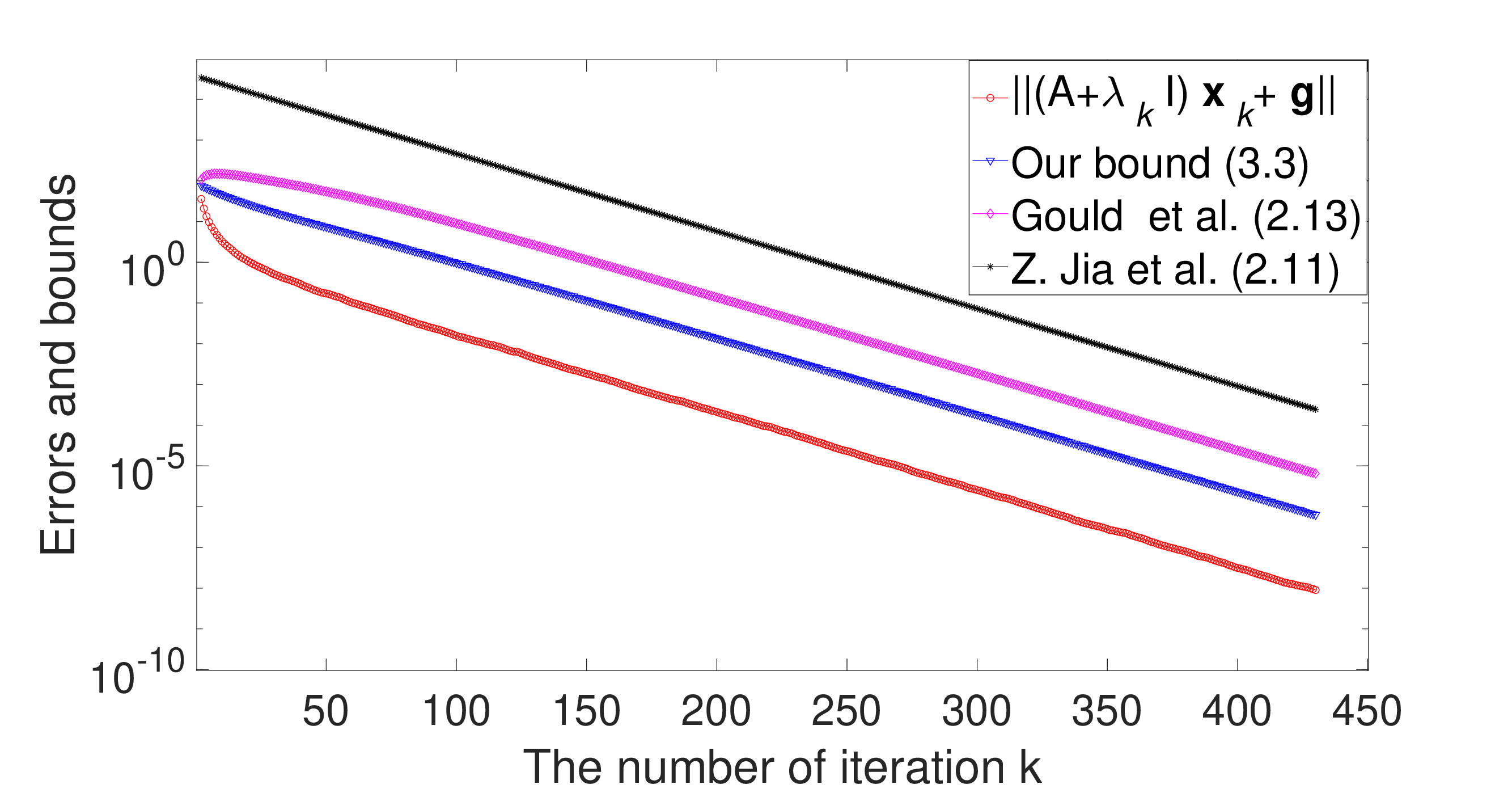}}
\\
\subfigure[\footnotesize  $A=diag(t^{[-50,50]}_{jn})$, $\D =15$.]{\includegraphics[width=6.3cm,height=3.8cm]{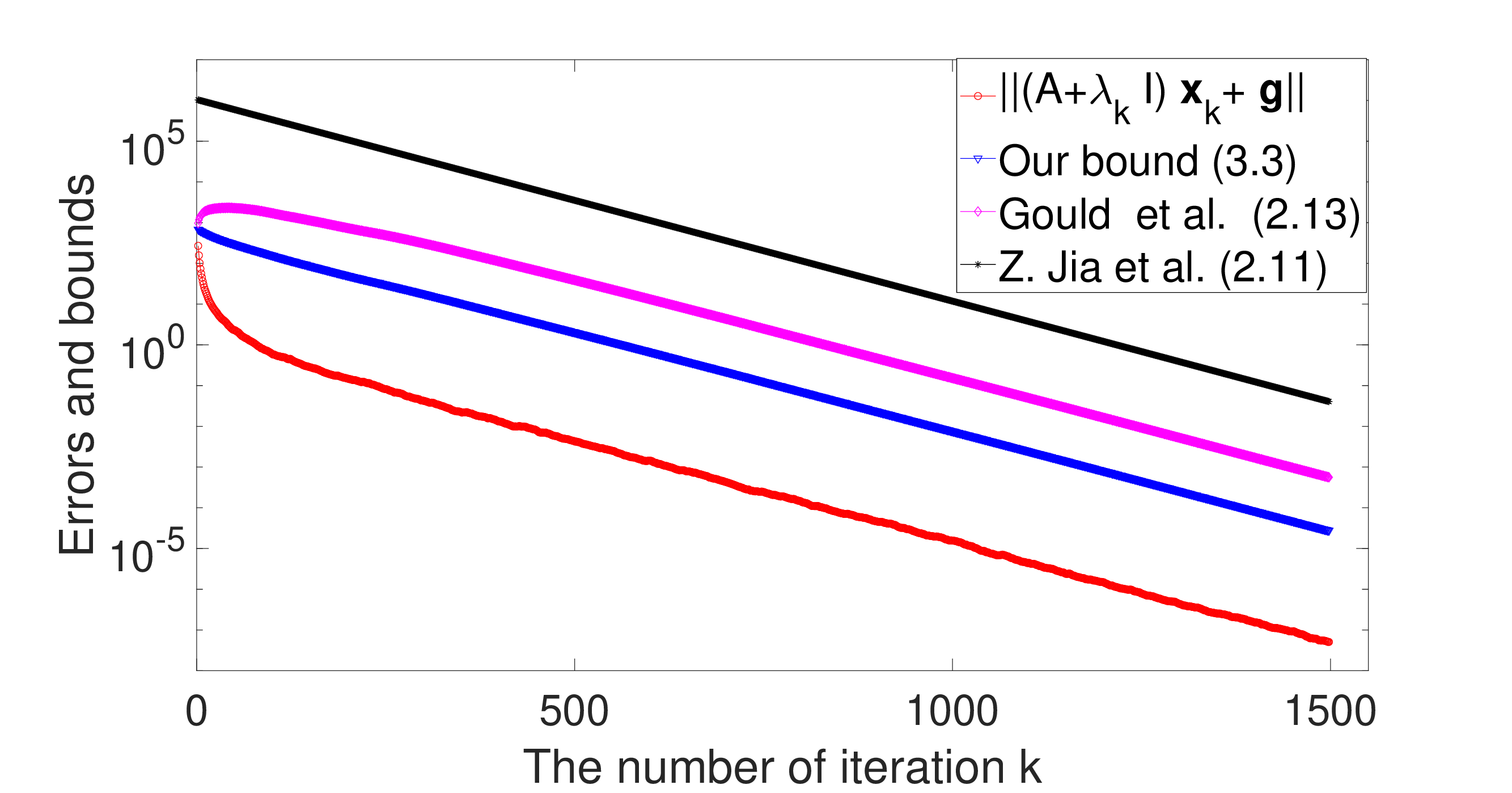}}~~~~
\subfigure[\footnotesize  $A=diag(t^{[-100,100]}_{jn})$, $\D
=20$.]{\includegraphics[width=6.3cm,height=3.8cm]{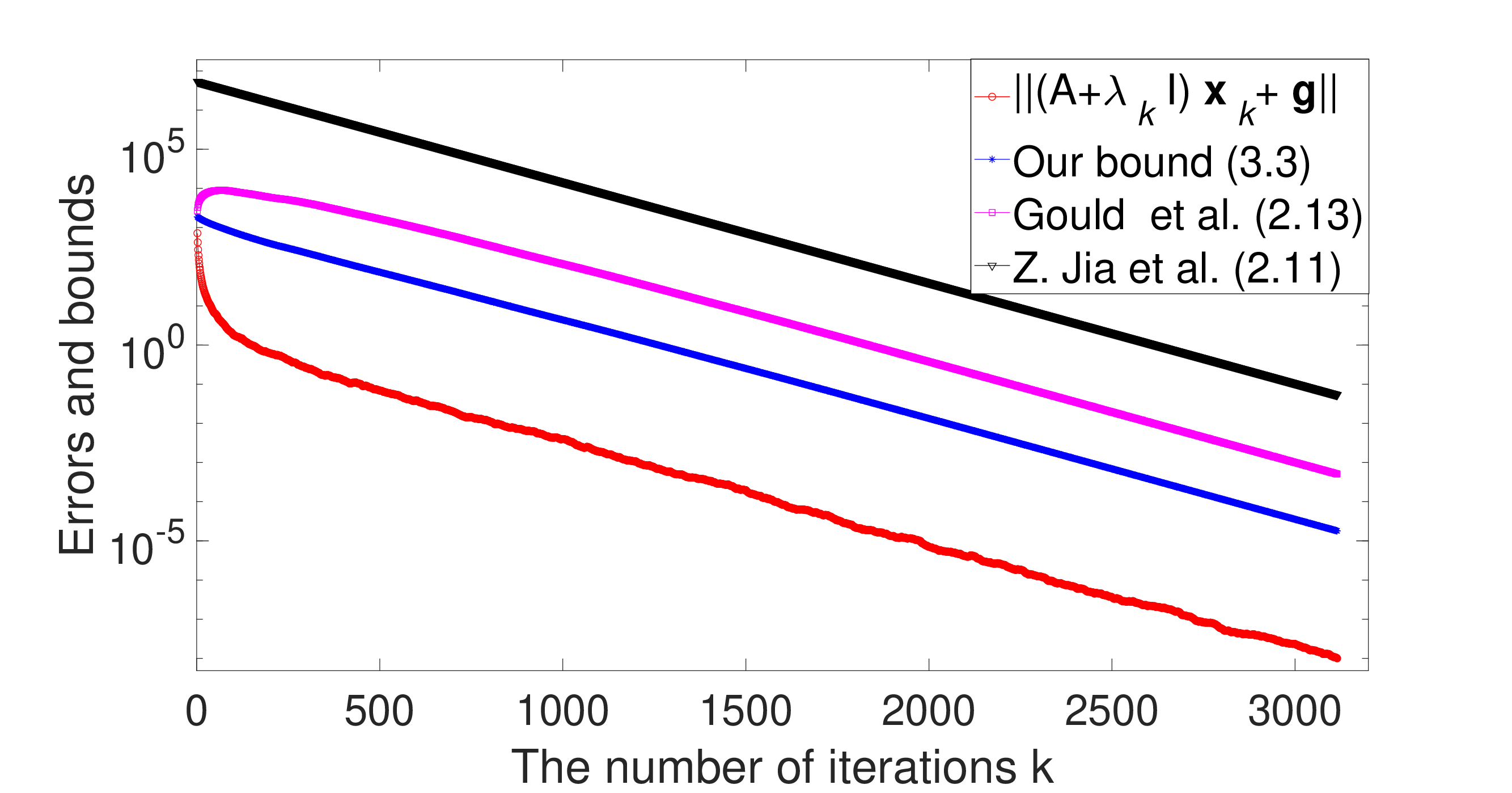}}
\end{figure}

\begin{figure}[ht]\caption{\footnotesize Example 1: A comparison of the upper bounds on  $\sin\angle(\bm x_{opt},\bm x_k)$.}\label{fig:1123}
\footnotesize{}
\centering
\subfigure[\footnotesize $A=diag(t^{[-5,5]}_{jn})$ with $\D
=1$.]{\includegraphics[width=6.3cm,height=3.8cm]{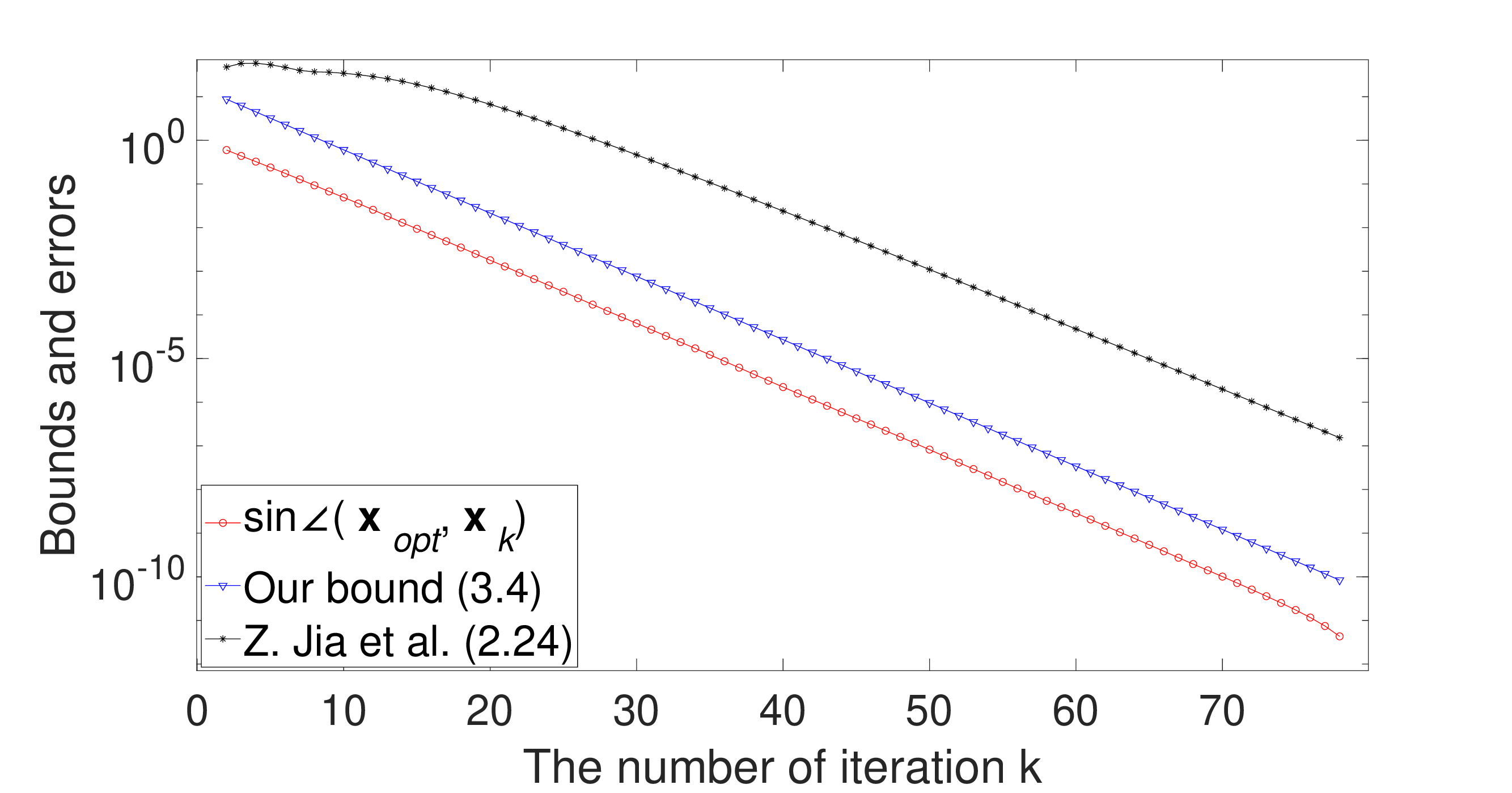}}~~~~
\subfigure[\footnotesize $A=diag(t^{[-10,10]}_{jn})$ with $\D
=10$.]{\includegraphics[width=6.3cm,height=3.8cm]{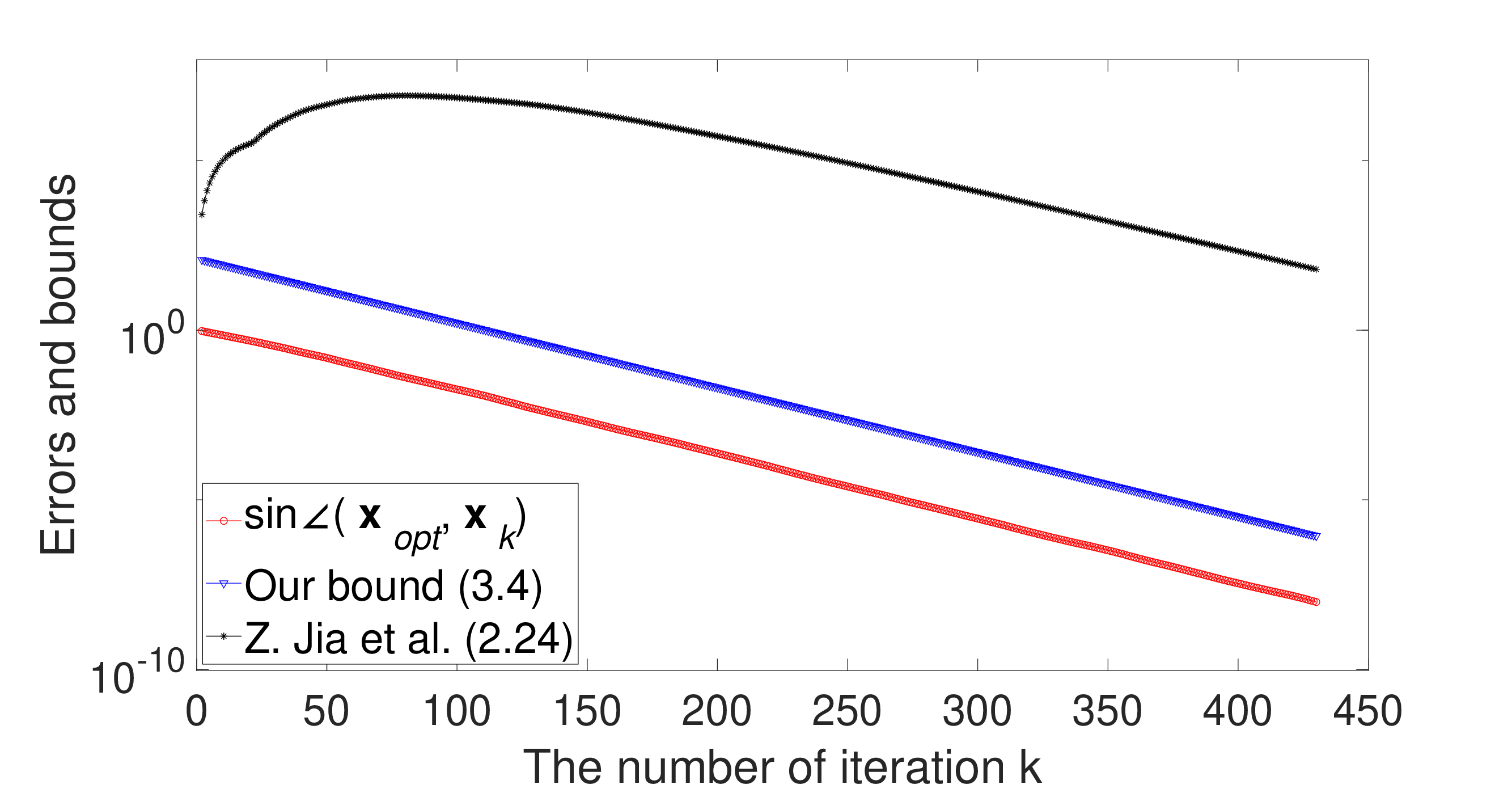}}
\\
\subfigure[\footnotesize $A=diag(t^{[-50,50]}_{jn})$ with $\D =15$.]{\includegraphics[width=6.3cm,height=3.8cm]{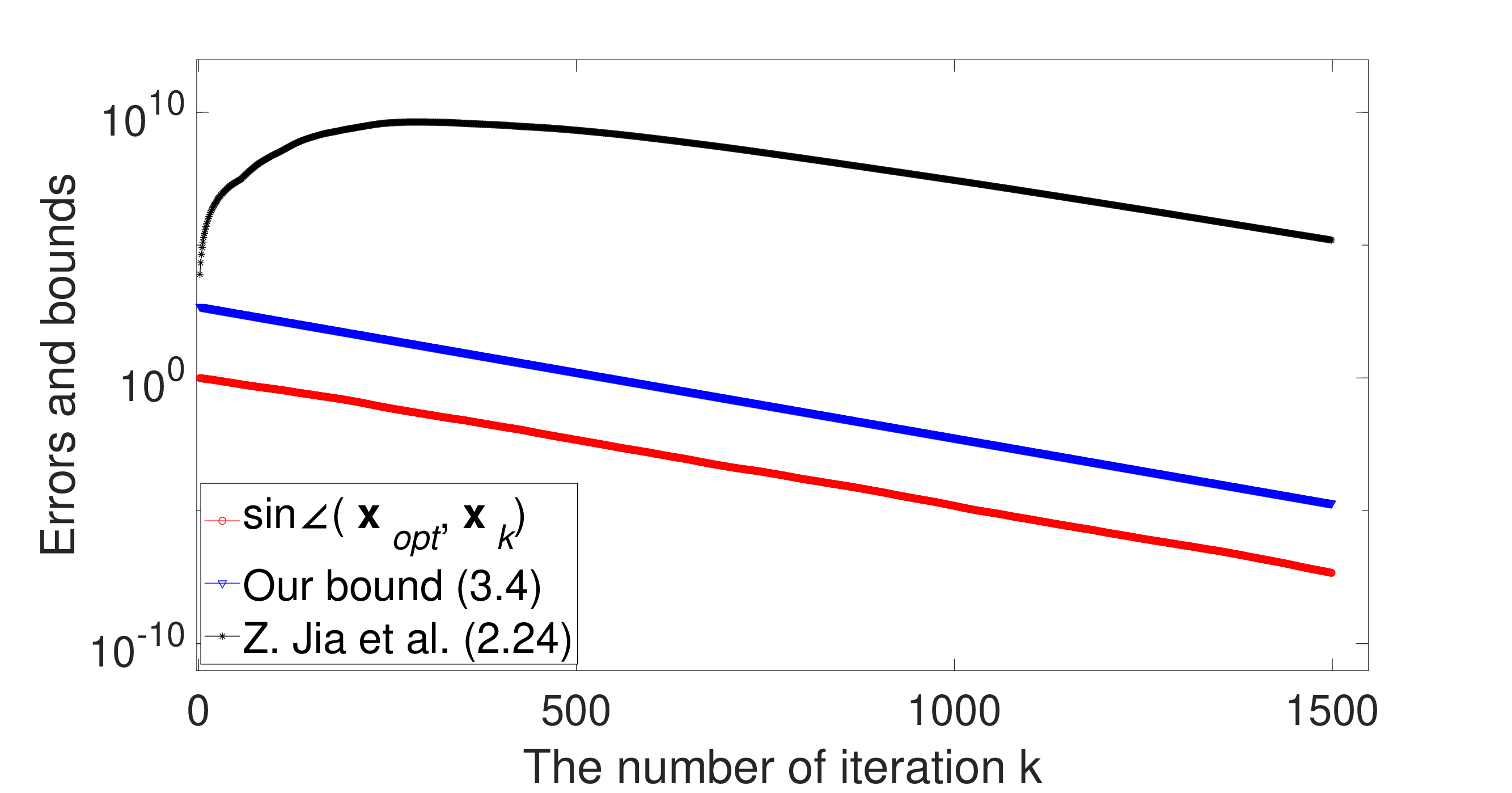}}~~~~
\subfigure[\footnotesize $A=diag(t^{[-100,100]}_{jn})$ with $\D
=20$.]{\includegraphics[width=6.3cm,height=3.8cm]{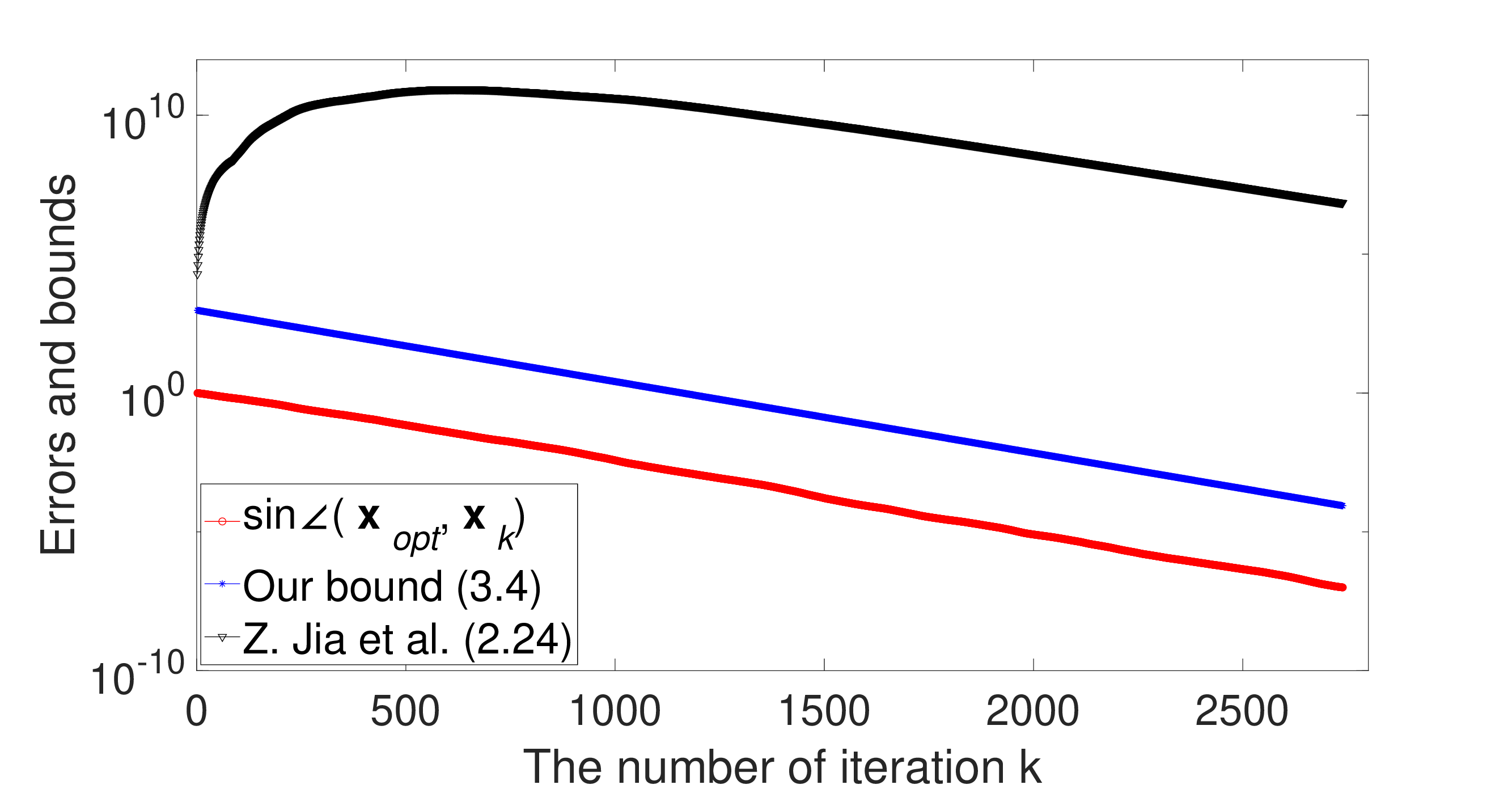}}
\end{figure}

\begin{figure}[ht]\caption{\footnotesize Example 1: A comparison of the upper bounds on  $\la_{opt}-\la_k$.}\label{19.1201}
\centering
\subfigure[\footnotesize $A=diag(t^{[-5,5]}_{jn})$ with $\D =1$.]{\includegraphics[width=6.3cm,height=3.8cm]{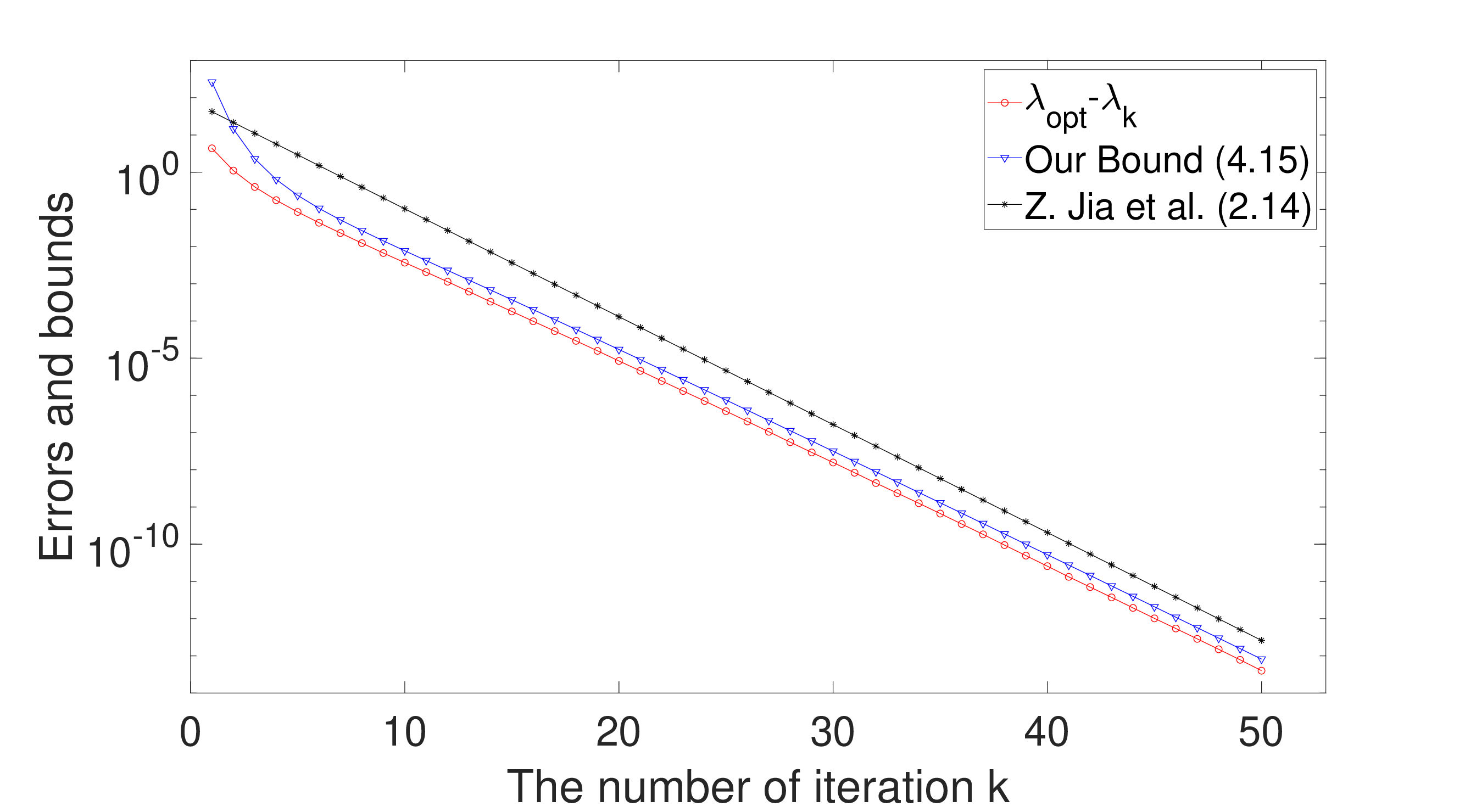}}~~~~
\subfigure[\footnotesize $A=diag(t^{[-10,10]}_{jn})$ with $\D =10$.]{\includegraphics[width=6.3cm,height=3.8cm]{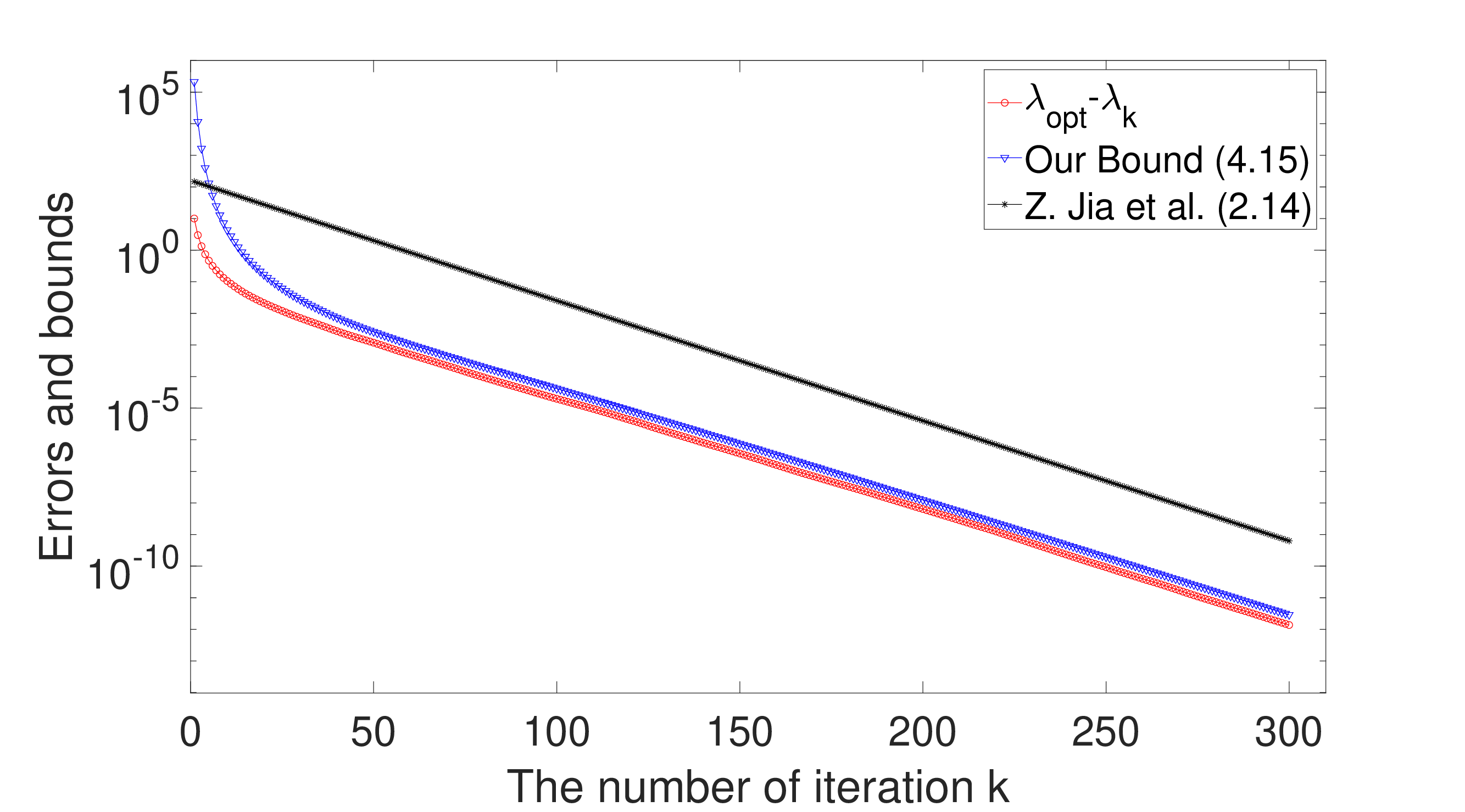}}
\\ 
\subfigure[\footnotesize $A=diag(t^{[-50,50]}_{jn})$ with $\D =15$.]{\includegraphics[width=6.3cm,height=3.8cm]{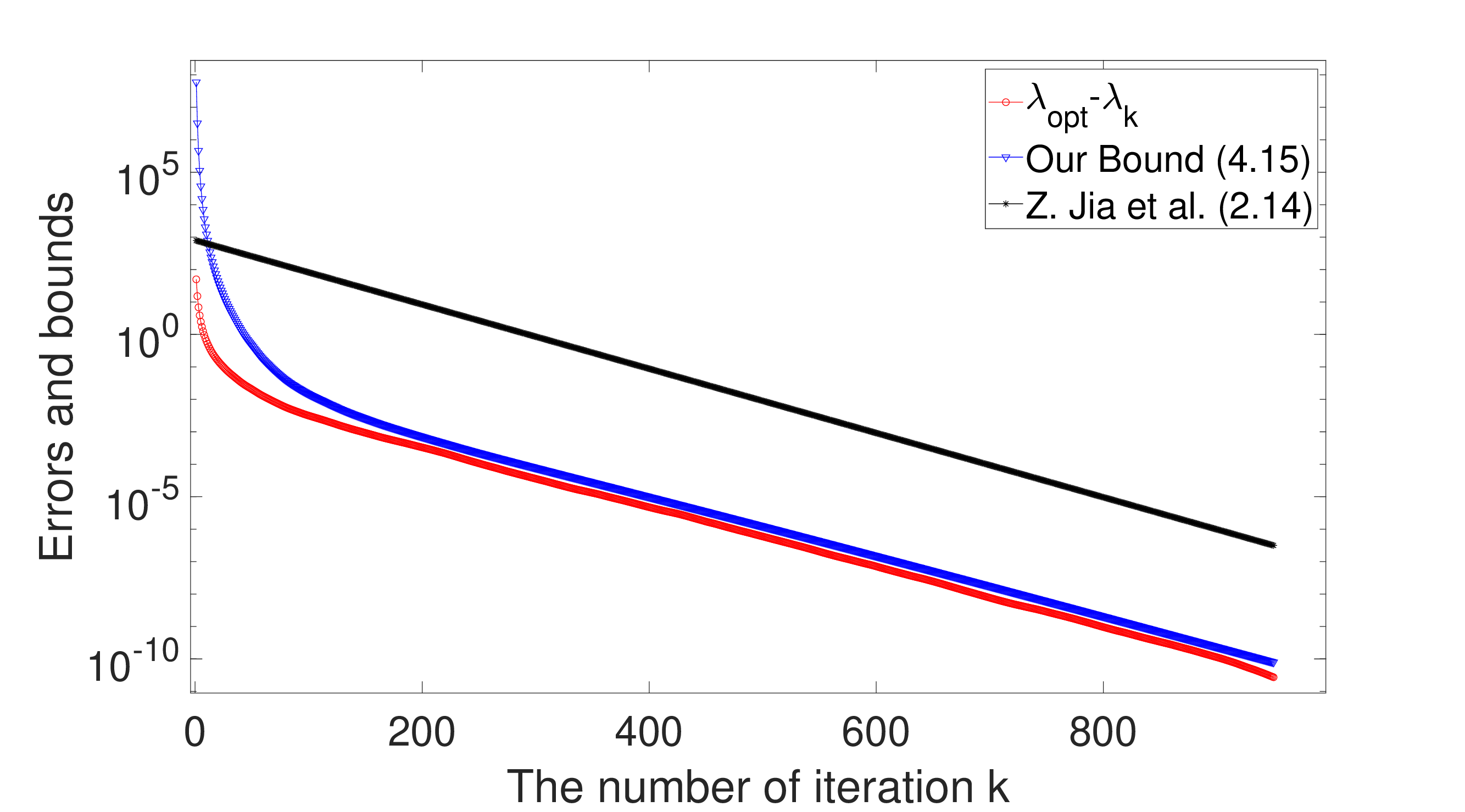}}~~~~
\subfigure[\footnotesize $A=diag(t^{[-100,100]}_{jn})$ with $\D =20$.]{\includegraphics[width=6.3cm,height=3.8cm]{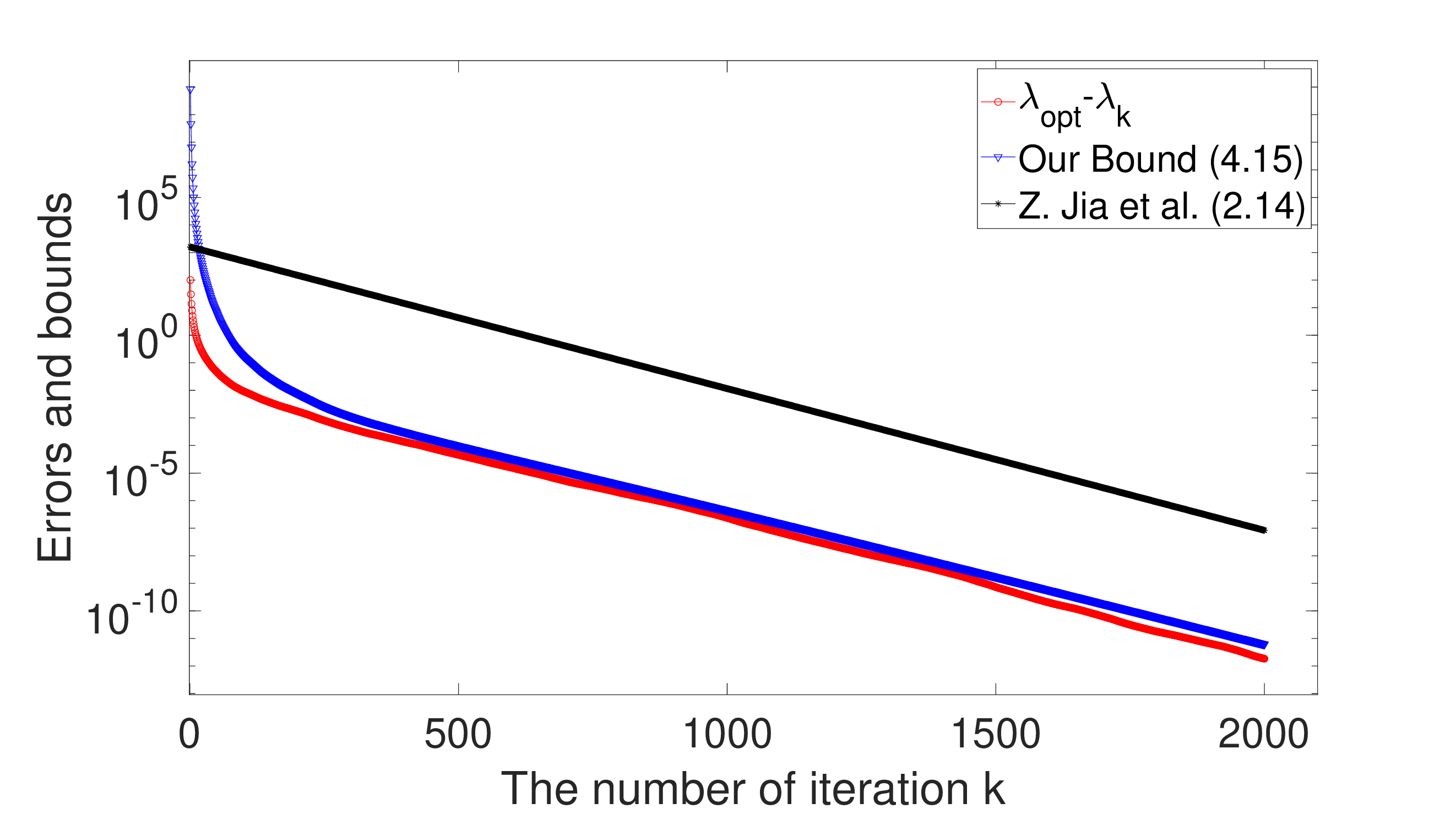}}
\end{figure}

In the upper bounds for comparison, the ``exact" values of $\la_{opt}$, $\bm x_{opt}$ and $\kappa$ are computed by using the $\tt GLTR$ method \cite{10}, which are listed in Table \ref{tab:13.54} with different $[a,b]$ and $\D$. Figure \ref{fig:1449}--Figure \ref{19.1201} plot the curves of $\|(A+\la_k I)\bm x_k+\bm g\|$, $\sin\angle(\bm x_{opt}, \bm x_{k})$, and $\la_{opt}-\la_k$, as well as the upper bounds for comparison as $k$ increases.

\begin{table}[ht]
\begin{center}\caption{\footnotesize Example 1: Approximations of $\la_{opt}$, $\|\bm x_{opt}\|$ and $\kappa$ compute by using $\tt GLTR$.}\label{tab:13.54}
\begin{tabular}{c|c|c|c|c|c}
\hline
  $[a,b]$&  $\D$  & $\la_{opt}$      &   $\|\bm x_{opt}\|$   &$\|(A+\la_{opt}I)\bm x_{opt}+\bm g\|$ & $\kappa$ \\
\hline
 [$-$5,~5]& $1$            &  5.2813 &     1.0000            &  $3.4226\times 10^{-12}$ & 36.5505  \\
\hline
 [$-$10,~10]  &  $10$       & 10.0096 &     10.0000       &   $9.0860\times 10^{-12}$ & $2.0887\times 10^{3}$\\
\hline
[$-$50,~50]    & 15         & 50.0032  &      15.0000       &  $9.7802\times 10^{-12}$ & $3.0794\times 10^{4}$ \\
\hline
[$-$100,~100]  &   $20$      &100.0018&  20.0000 &   $3.9402\times 10^{-11}$ & $1.1428\times 10^{5}$ \\
\hline
\end{tabular}
\end{center}
\end{table}
\begin{table}[ht]
\begin{center}\caption{\footnotesize Example 1: A comparison of $s(\la_{opt})=\frac{\D^2}{\bm g^TA_{opt}^{-3}\bm g}$  and $cond(\la_{opt})=\frac{1}{2}\vert\bm y_1^T\bm y_2\vert^{-1}$.}\label{tab:16.39}
\begin{tabular}{c|c|c|c}
\hline
  $[a,b]$&  $\D$  & $s(\la_{opt})$      &   $cond(\la_{opt})$   \\
\hline
 [$-$5,~5]& $1$            &  0.3811 &     1.6073           \\
\hline
 [$-$10,~10]  &  $10$       & 0.0129 &     43.6256    \\
\hline
[$-$50,~50]    & 15         & 0.0044  &      1.2807$\times 10^2$      \\
\hline
[$-$100,~100]  &   $20$      &0.0023&  2.3780$\times 10^2$  \\
\hline
\end{tabular}
\end{center}
\end{table}

It is obvious to see from the figures that our new upper bounds are much sharper than those of Z. Jia {\it et al}. \cite{5} and  Gould {\it et al.} \cite{G.S.}.
First, it is seen from Figure \ref{fig:1449} that for the upper bounds $\|(A+\la_k I)\bm x_k+\bm g\|$, our bound
\eqref{10.571} is about 10 times smaller than \eqref{10.28} (the one due to Gould {\it et al.}), and is about 1000 times smaller than \eqref{10.27} (the one due to Z. Jia {\it et al}.), especially when the condition number $\kappa$ is large. 
Second,
Figure \ref{fig:1123} demonstrate that the new upper bound of $\sin\angle(\bm x_{opt},\bm x_k)$ is much smaller than the one due to Z. Jia {\it et al}. This is because $sep(\la_{opt},C_k)$
is usually very small in practice, and the new upper bound \eqref{16.0600} can be much smaller than \eqref{11.2300}; refer to \eqref{eq47}.
Furthermore, we observe from Figure \ref{fig:1123} that the curve of the upper bound established by Z. Jia {\it et al}. is first up and then down as $k$ increases, especially when $\kappa$ is relatively large. This is because one has to first offset the value of $\widetilde{c}_k$ (which is also increased on $k$, see \eqref{eq47}), by using $\big(\frac{\sqrt{\kappa}-1}{\sqrt{\kappa}+1}\big )^k$.

Third, we see from Figure \ref{19.1201} that for $\la_{opt}-\la_k$, our bound \eqref{8.26} is much better than \eqref{11.15} due to Z. Jia {\it et al}.
However, it seems that the new upper bound is a little large at the beginning of the iterations. This is due to the fact that $\|\bm g\|^2\cdot\bm e_1^T(T_k+\la_{opt} I)^{-3}\bm e_1$, the denominator of $s_k$ (see \eqref{eq513}), can be small at the beginning. Note that this coincides with the trend of the convergence of $\la_{opt}-\la_k$.

Furthermore, it is shown from Theorem \ref{8.42} that $s(\la_{opt})$ plays an important role in the convergence of $\lambda_k$.
In Remark \ref{remk16.59}, we pointed out $s(\la_{opt})$ is closely related to the condition number $cond(\la_{opt})$, however,
the former can be much smaller than that of the latter. In order to show this more precisely, we list in Table \ref{tab:16.39} the two values for different
$[a,b]$ and $\Delta$. One observes that $s(\la_{opt})$ can be about $10^4$ times smaller than $cond(\la_{opt})$ in practice. One refers to Example 2 for more details on the importance of $s(\la_{opt})$.

{\bf Example 2.}
In this example, we try to show the importance of
$$
s(\la_{opt})=\frac{\D^2}{\bm g^TA_{opt}^{-3}\bm g}
$$
on the convergence of $\la_k$; refer to \eqref{eq513}.
Consider the matrix
$$
A_1=diag(t^{[a,b]}_{jn}),~j=1,2,\ldots,n,~~{\rm with}~a=1,~b=3000,~{\rm and}~n=10000,
$$
where $t^{[a,b]}_{jn}$ is defined in \eqref{eqn61}.
Let $\mathcal{A}_1=diag(A_1,0.01)$, and
$$
{A}=\mathcal{A}_1-500 I,
~~{\rm and}~\bm f_\zeta=(\zeta,\ldots,\zeta,1)^T\in \mathbb{R}^{10001},
$$
where $I$ is identity matrix of appropriate size.
We set
$$
{\D=\|\mathcal{A}_1\bm f_\zeta\|}~~{\rm and}~~\bm g=\mathcal{A}_1^2\bm f_\zeta.
$$
Thus, for any $\zeta\in\mathbb{R}$, we have that
$$
\la_{opt}=500,~\bm x_{opt} = -\mathcal{A}_1\bm f_\zeta,~A_{opt}=\mathcal{A}_1~{\rm and}~\kappa(A_{opt})=3\times 10^5.
$$

\begin{figure}[ht]
    \centering
	\begin{minipage}{0.45\linewidth}
		\centering
		\vspace{0cm}
		\setlength{\abovecaptionskip}{0.28cm}
		\includegraphics[width=6.3cm,height=3.8cm]{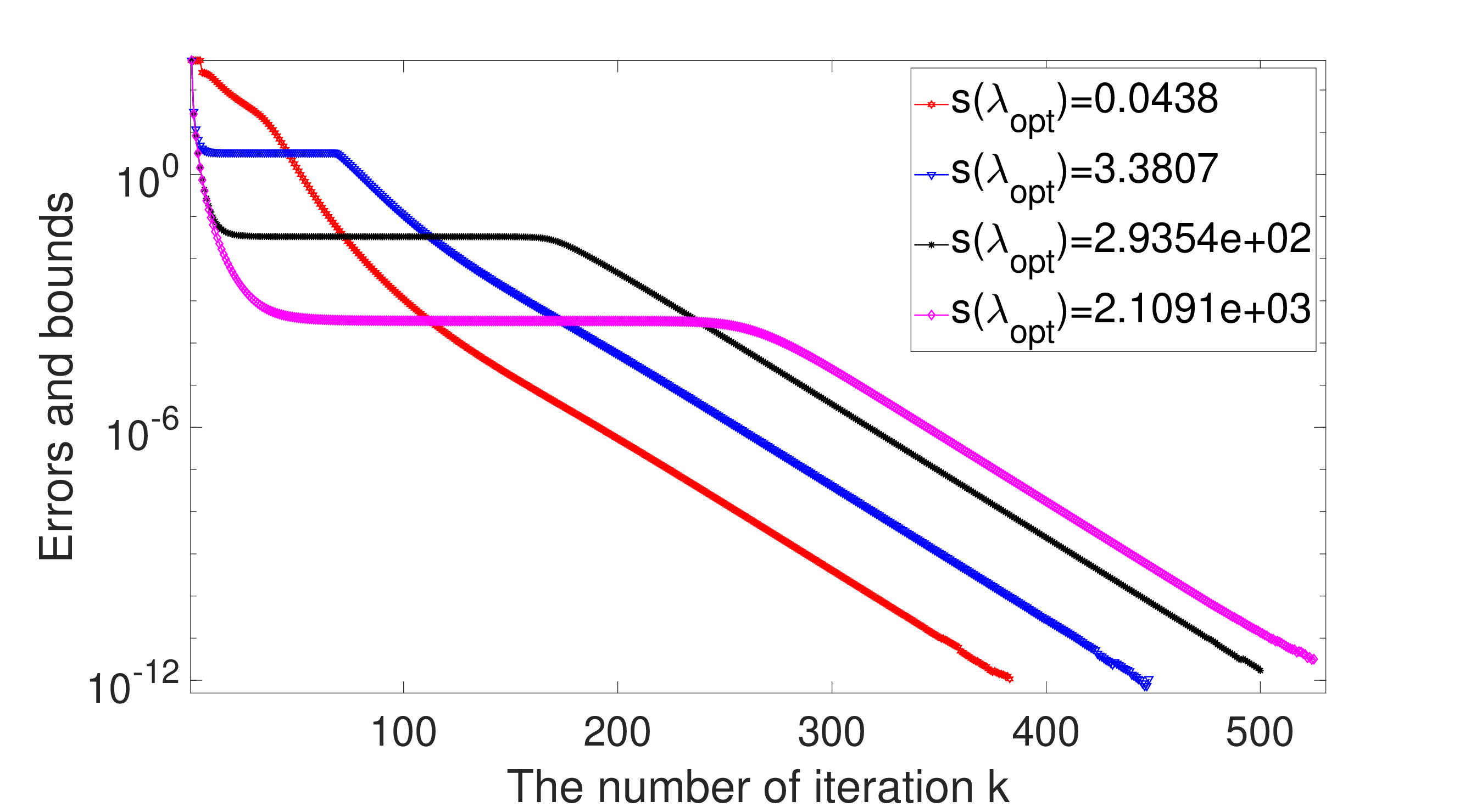}
		\caption{\footnotesize Example 2: Curves of $\la_{opt}-\la_{k}$ with different $s(\la_{opt})$.}
		\label{fig19.1201}
	\end{minipage}
	\hfill
	\begin{minipage}{0.45\linewidth}
		\centering
        \begin{tabular}{c|c}
\hline
$\zeta$& $s(\la_{opt})$\\
\hline
$10^{-4}$&$2.1091\times 10^{3}$\\
\hline
$10^{-5}$&$2.9354e\times 10^{2}$\\
\hline
$10^{-6}$& 3.3807\\
\hline
$10^{-7}$& 0.0438\\
\hline
\end{tabular}
\captionof{table}{\footnotesize Example 2: The values of $s(\la_{opt})$ with different $\zeta$.}\label{tab:21.3444}
	    \end{minipage}
\end{figure}

To illustrate the importance of ${s}(\la_{opt})$ on the convergence of $\la_k$, we present in Table \ref{tab:21.3444} the values of ${s}(\la_{opt})$ corresponding to different $\zeta$. In Figure \ref{fig19.1201}, we plot the curves of $\la_{opt}-\la_k$ with different ${s}(\la_{opt})$. It is obvious to see from the figure that the larger ${s}(\la_{opt})$, the slower $\la_k$ converges. Recall from Theorem \ref{8.42} and Theorem \ref{thm20.08} that ${s}(\la_{opt})$ and $\kappa$ are two key factors for the convergence of $\la_k$. As the condition numbers $\kappa$ are the same in all the cases, the differences are mainly from those of ${s}(\la_{opt})$. This shows the effectiveness of our theoretical results, refer to Remark \ref{14.55}.

{\bf Example 3.}
In this example, we consider the cubic regularization problem \eqref{2}.
To show the effectiveness of our theoretical results, we compare our upper bound \eqref{eq.8.26} with \eqref{eq2054} established by X. Jia \emph{et al.}
Consider the matrix used in Example 1, where
$$
A=diag(t^{[a,b]}_{jn}),~j=1,2,\ldots,n,~~{\rm with}~a=-1,~b=1,~{\rm and}~n=10000.
$$

\begin{table}[ht]
\begin{center}\caption{\footnotesize Example 3: $\sigma$, $\mu_{opt}$, $\kappa_{\sg,k}$ and residual norms in the cubic regularization problem.}\label{tab:21360}
\begin{tabular}{c|c|c|c}
\hline
 $\sigma$    &  $\mu_{opt}$  &  $\kappa_{\sg,k}$ & ${\|(A+\mu_{opt}I)\bm x_{opt}+\bm g\|}$\\
\hline
 1  &          1.3405&   6.8729& 1.3439e-16   \\
\hline
0.1&             1.0214 &94.5315   &  7.3994e-16 \\
 \hline
    0.01&  1.0010  &  2.0914e+03 & 6.1067e-15  \\
\hline
  0.001        & 1.0000 &   4.5269e+04  &  7.9361e-14 \\
   \hline
\end{tabular}
\end{center}
\end{table}
It was pointed out that \eqref{2} can be equivalently rewritten as a large-scale eigenvalue problem \cite{F.L}.
In this example, we make use of the MATLAB built-in function $\tt eigs.m$ (with stopping tolerance $tol=10^{-12}$) to compute $\mu_{opt}$ and $\mu_k$. As $A_{opt}$ is a diagonal matrix in this example, the solution $\bm x_{opt}^\sg$ is computed by using the dot division commend ``$.\setminus$" in MATLAB.

\begin{figure}[ht]\caption{\footnotesize Example 3: A comparison of the upper bounds on $\mu^3_{opt}-\mu_k^3$.}\label{fig:2124}
\centering
\subfigure[\footnotesize $\sigma=1$.]{\includegraphics[width=6.3cm,height=3.8cm]{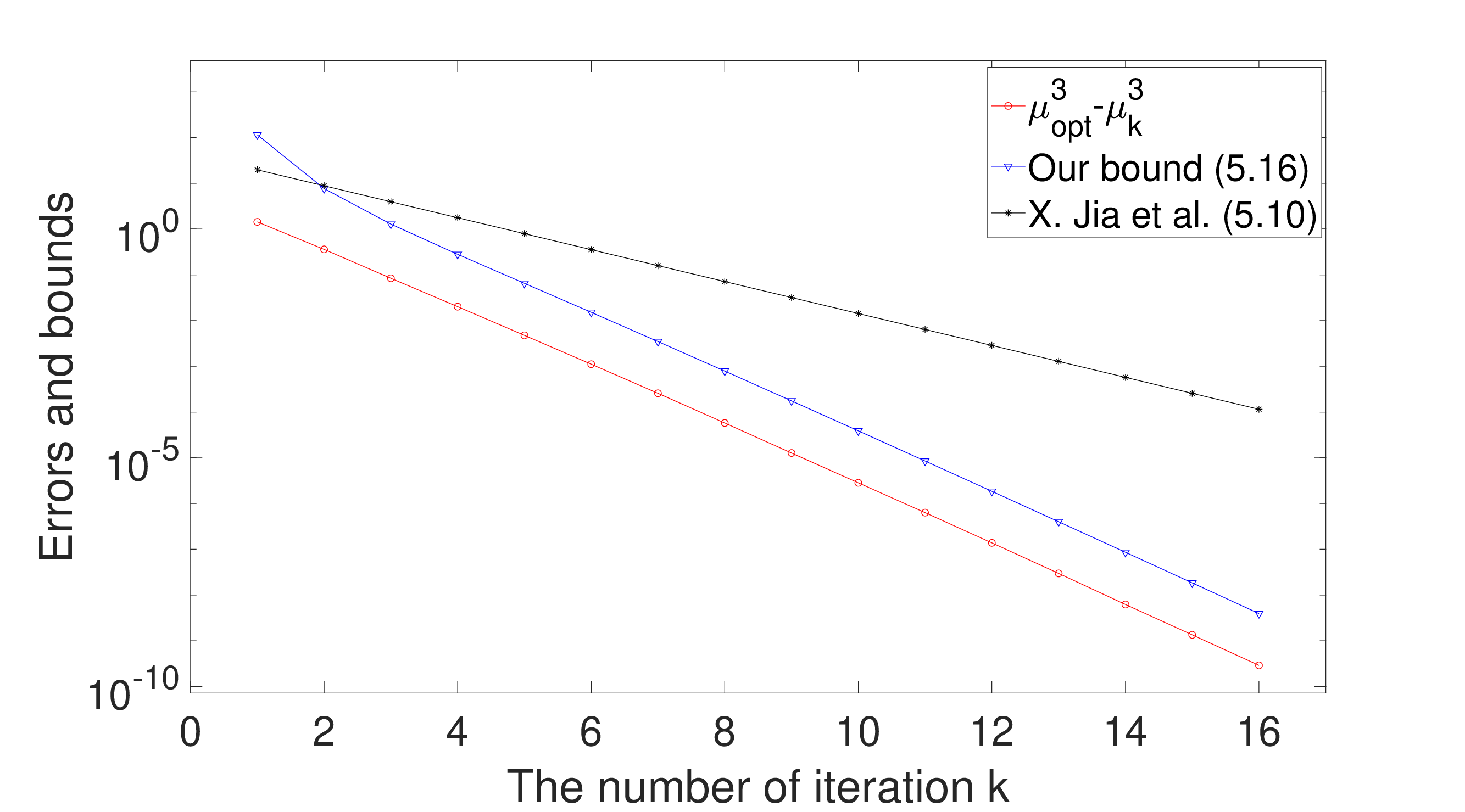}}~~~~
\subfigure[\footnotesize $\sigma=0.1$.]{\includegraphics[width=6.3cm,height=3.8cm]{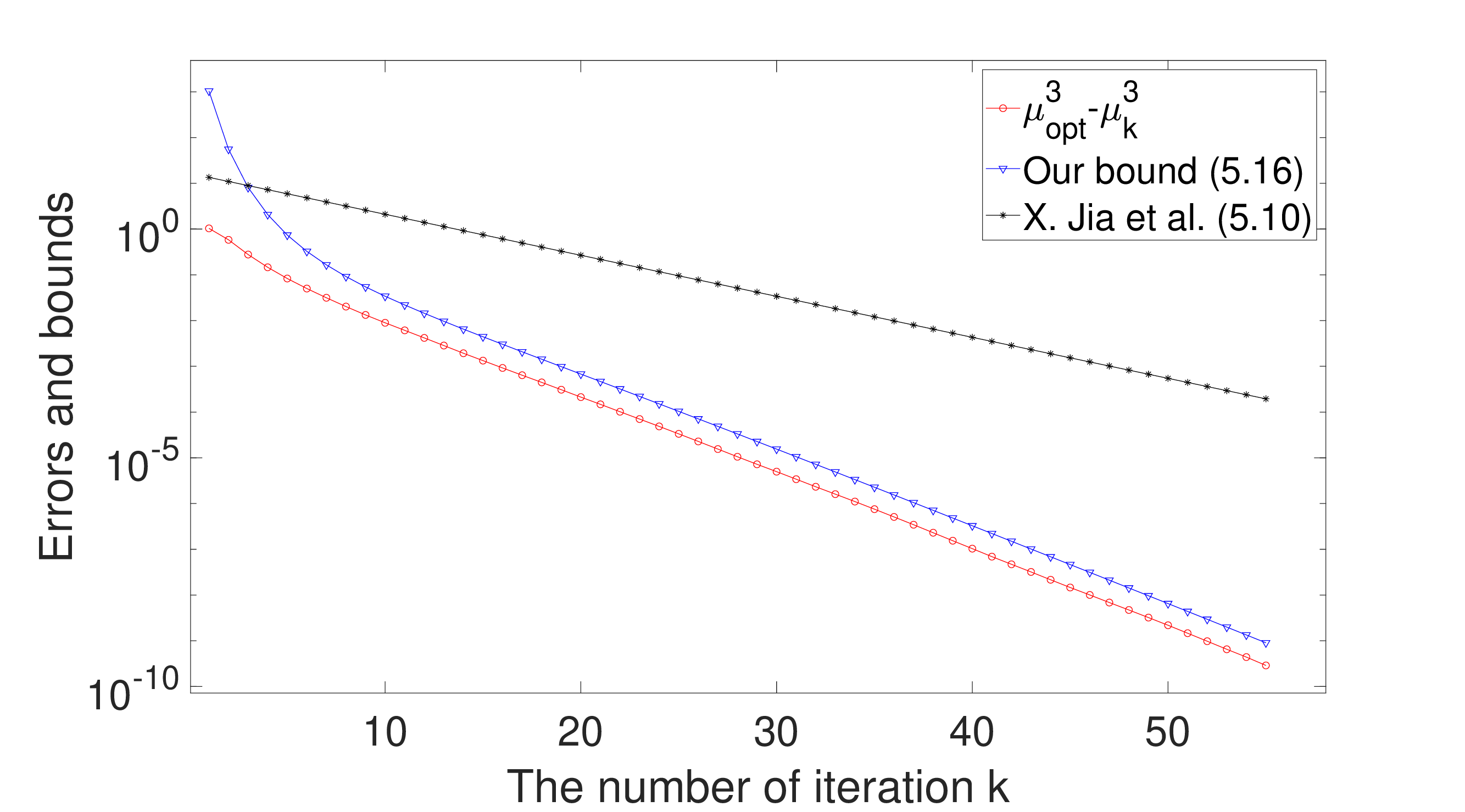}}
\\
\subfigure[\footnotesize $\sigma=0.01$.]{\includegraphics[width=6.3cm,height=3.8cm]{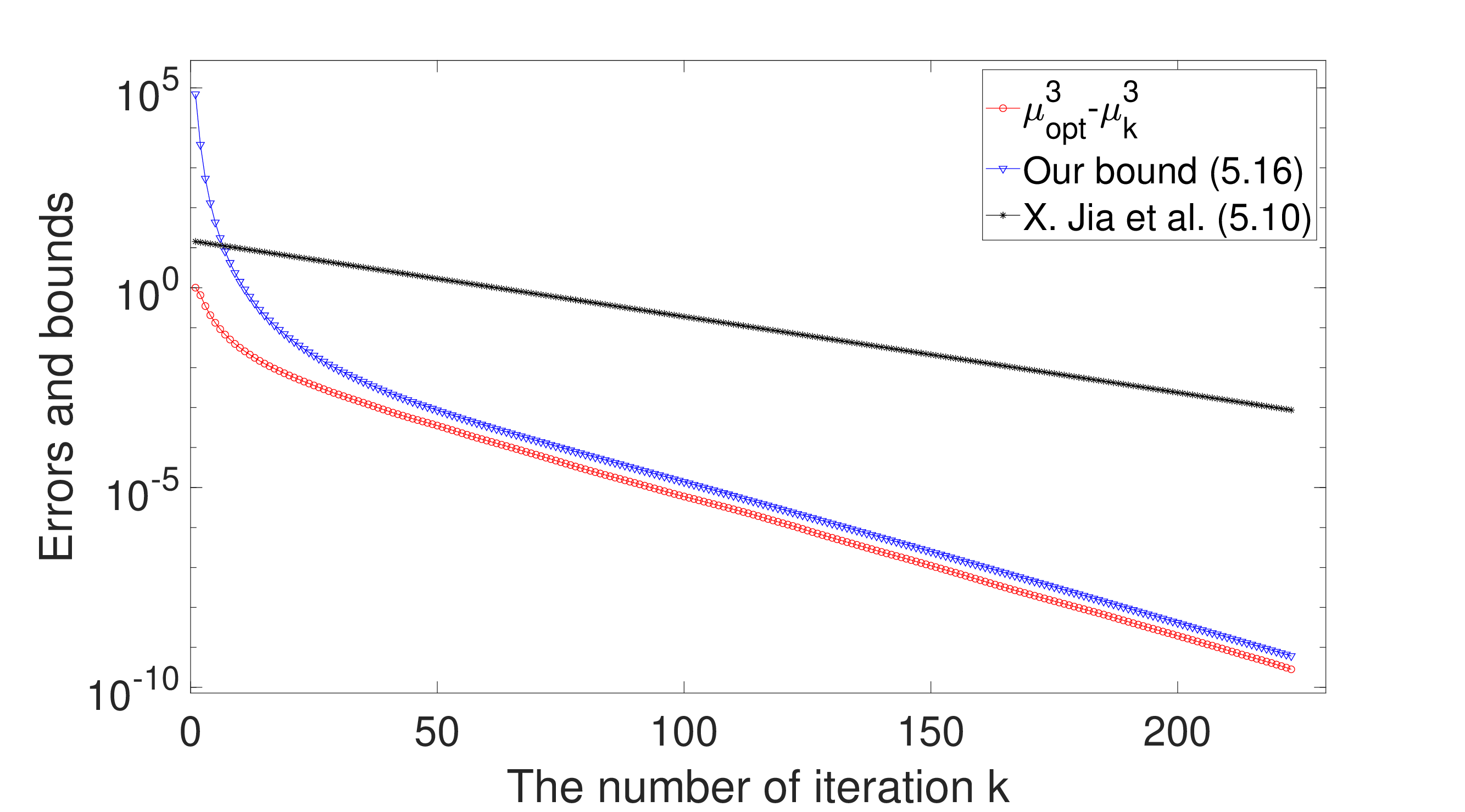}}~~~~
\subfigure[\footnotesize $\sigma=0.001$.]{\includegraphics[width=6.3cm,height=3.8cm]{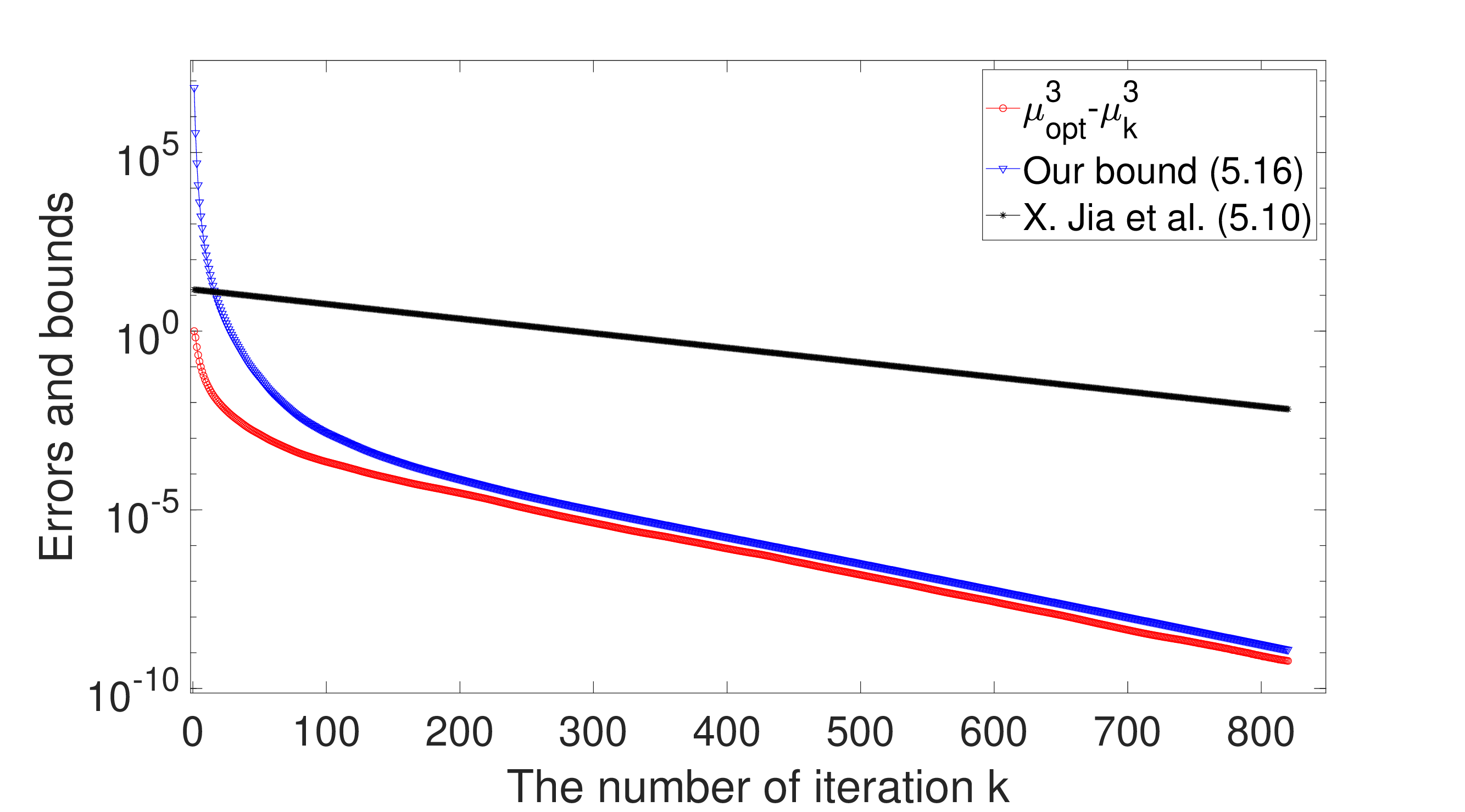}}
\end{figure}

The choices of $\sigma$, the values of $\mu_{opt}$, $\kappa_{\sg,k}$ and residual norms are listed in Table \ref{tab:21360}.
Figure \ref{fig:2124} plots the curves of $\mu^3_{opt}-\mu_k^3$ and the two upper bounds as $k$ increases. Two remarks are in order. First,
we see from Figure \ref{fig:2124} that our result is much sharper than that of X. Jia \emph{et al.}, and the convergence rates of the two upper bounds are different in essence. More precisely, the convergence rate of our new bound  is $\big(\frac{\sqrt{\kappa_{\sigma}}-1}{\sqrt{\kappa_{\sigma}}+1}\big )^{2}$, while that of X. Jia \emph{et al.} is $\big(\frac{\sqrt{\kappa_{\sigma}}-1}{\sqrt{\kappa_{\sigma}}+1}\big )$; please see \eqref{eq.8.26} and \eqref{eq2054}, respectively.
Second,
similar to Figure \ref{19.1201}, it is seen from Figure \ref{fig:2124} that the new upper bound is relatively large at the initial stage of the iteration, which coincides with the trend of the real values of $\mu^3_{opt}-\mu^3_k$.
This is because the value of $\|\bm g\|^2\cdot\bm e_1^T(T_k+\mu_{opt} I)^{-3}\bm e_1$ can be small in the initial stage.

\section{Conclusion}
The GLTR method is a popular approach for solving large-scale TRS. In essence, this method is a
projection method in which the original large-scale TRS is projected into
a small-sized one. Recently, Z. Jia {\it et al.} considered the convergence of the GLTR method \cite{5}.
In this paper, we revisit this problem and establish some refined bounds on the residual norm $\|(A+\la_k I)\bm x_k+\bm g\|$, the distance $\sin\angle(\bm x_{opt}, \bm x_k)$ between the exact solution $\bm x_{opt}$ and the approximate solution $\bm x_k$, as well as the distance $\la_{opt}-\la_k$ between the Lagrange multiplier $\la_{opt}$ and its approximation $\la_k$. Moreover, we generalize these results to the convergence of Krylov subspace method for the cubic regularization problem, and improve some bounds due to X. Jia {\it et al} \cite{J&Z}.

In this paper, we only consider the situation of {\it easy case} of the TRS \eqref{1}. Further study includes the perturbation analysis on this problem, and proposing more efficient methods for judging and solving the TRS \eqref{1} in the ({\it nearly}) {\it hard} case. These are very interesting topics and are definitely a part of our future work.

\end{document}